RICE UNIVERSITY

# Spectral Properties of Limit-Periodic Schrödinger Operators

by

**Zheng Gan**

A Thesis Submitted
in Partial Fulfillment of the
Requirements for the Degree

**Doctor of Philosophy**

Approved, Thesis Committee:

David Damanik,
Professor of Mathematics, Chair

Robert Hardt,
W. L. Moody Professor of Mathematics

Mark Embree,
Professor of
Computational and Applied Mathematics

Houston, Texas
Nov-28-2011

# Abstract

Spectral Properties of Limit-Periodic Schrödinger Operators
by
Zheng Gan

We investigate spectral properties of limit-periodic Schrödinger operators in $\ell^2(\mathbb{Z})$. Our goal is to exhibit as rich a spectral picture as possible. We regard limit-periodic potentials as generated by continuous sampling along the orbits of a minimal translation of a procyclic group. This perspective was first proposed by Avila and further exploited by the author, which allows one to separate the base dynamics and the sampling function. Starting from this point of view, we conclude that all the spectral types (i.e. purely absolutely continuous, purely singular continuous, and pure point) can appear within the class of limit-periodic Schrödinger operators. We furthermore answer questions regarding how often a certain type of spectrum would occur and discuss the corresponding Lyapunov exponent. In the regime of pure point spectrum, we exhibit the first almost periodic examples that are uniformly localized across the hull and the spectrum.

## Acknowledgments

I am deeply grateful to my advisor, David Damanik, whose encouragement, guidance and support enabled me to develop an understanding of the subject and finally finish the thesis project.

I am grateful to my academic brother, Helge Krüger, for never ending mathematical discussions.

I am grateful to many professors in our Schrödinger operator research community, like Gunter Stolz, Svetlana Jitomirskaya, Christian Remling, Yoram Last, Barry Simon, etc., who always gave me suggestions and answered my questions whenever I met them in conferences.

Last, I wish to thank all of the professors and staffs at the math department. They work hard for the department, creating a productive and pleasant studying environment. I wish to thank all of the fellow graduate students, who make our department as a warm and fun family.

To my parents.

# Contents

# Introduction

In this thesis, we will study spectral properties of limit-periodic Schrödinger operators. Someone may ask, why do people care about spectral properties of Schrödinger operators? Let's start with the celebrated Schrödinger equation, which is used to describe how the quantum state of a physical system evolves in time.

In quantum mechanics, a particle moving around in a $d$-dimensional space is described by a normalized vector $\psi(t)$ in the separable Hilbert space $\mathcal{H} = L^2(\mathbb{R}^d)$ ($\psi(t)$ is also called a wavefunction or state) which satisfies the equation

$$i\frac{\partial}{\partial t}\psi = H\psi, \tag{1}$$

where $H$, called the Hamiltonian, is a self-adjoint operator on $L^2(\mathbb{R}^d)$. In general, we assume that $H$ is independent of time $t$, which is realizable in most physical models. We can solve the above equation as

$$\psi(t) = e^{-itH}\psi,$$

where $\psi$ is the state of the quantum particle at $t = 0$. Then, a natural question arises: how does one describe $\psi(t) = e^{-itH}\psi$?

Assume that $\{\delta_n\}_{n\in\mathbb{Z}}$ is an orthonormal basis of $\mathcal{H}$. Then $|\langle\delta_n, \psi(t)\rangle|^2$ can be interpreted as the probability of finding the particle at the state $\delta_n$ at time $t$. The famous spectral theorem is the following.

**Theorem** (Spectral Theorem)**.** *For every bounded measurable function $f : \mathbb{R} \to \mathbb{R}$, there exists a unique measure $\mu_\psi$ such that*

$$\langle\psi, f(H)\psi\rangle = \int_{\sigma(H)} f(x)\ d\mu_\psi.$$

There are many variants of the spectral theorem. Please refer to [45] for more details about the spectral theorem. With the spectral theorem, we have

$$\mathcal{H}_{ac} = \{\psi \in \mathcal{H} \mid \mu_\psi \text{ is absolutely continuous with respect to the Lebesgue measure}\},$$

$$\mathcal{H}_{sc} = \{\psi \in \mathcal{H} \mid \mu_\psi \text{ is singular continuous with respect to the Lebesgue measure}\},$$

and

$$\mathcal{H}_{pp} = \{\varphi \in \mathcal{H} \mid \mu_\psi \text{ is pure point with respect to the Lebesgue measure}\}.$$

It is easy to conclude $\mathcal{H} = \mathcal{H}_{ac} \oplus \mathcal{H}_{sc} \oplus \mathcal{H}_{pp}$. Let $\sigma_{ac}(H) = \sigma(H|_{\mathcal{H}_{ac}})$, $\sigma_{sc}(H) = \sigma(H|_{\mathcal{H}_{sc}})$, and $\sigma_{pp}(H) = \sigma(H|_{\mathcal{H}_{pp}})$. We have $\sigma(H) = \sigma_{ac}(H) \cup \sigma_{sc}(H) \cup \sigma_{pp}(H)$ (note that they don't need to be pairwise disjoint).

**Theorem** (RAGE Theorem, Ruelle-Amrein-Georgescu-Enss)**.** *Assume $\psi \in \mathcal{H}$ and $\psi(t)$ is the solution of the Schrödinger equation (1). Then we have*
*(1). $\mu_\psi$ is pure point if and only if for any $\varepsilon > 0$, there is $N \in \mathbb{Z}^+$ such that*

$$\sum_{|n| \geq N} | \langle \delta_n, \psi(t) \rangle |^2 < \varepsilon$$

*for every $t \in \mathbb{R}$.*
*(2). $\mu_\psi$ is purely continuous if and only if for every $N \in \mathbb{Z}^+$,*

$$\lim_{T \to \infty} \frac{1}{T} \int_{-T}^{T} \sum_{|n| \leq N} | \langle \delta_n, \psi(t) \rangle |^2 \, dt = 0.$$

*(3). If $\mu_\psi$ is purely absolutely continuous, then for every $N \in \mathbb{Z}^+$,*

$$\lim_{|t| \to \infty} \sum_{|n| \leq N} | \langle \delta_n, \psi(t) \rangle |^2 = 0.$$

The RAGE theorem asserts, roughly speaking, that time evolution of $\psi(t)$ is closely related to the spectral types of $H$. In $\mathcal{H} = L^2(\mathbb{R}^d)$, for a free particle,

$$H = -\Delta = -\sum_{j=1}^{d} \frac{\partial^2}{\partial x_j^2},$$

which presents the kinetic energy. If there is a potential energy $V(x)$, we may also use

$$(H\psi)(x) = -\Delta\psi(x) + V(x)\psi(x).$$

We call the $H$ a $d$-dimensional continuum Schrödinger operator.

If $\mathcal{H} = \ell^2(\mathbb{Z}^d)$, then $-\Delta$ is replaced by a finite difference operator and $V(x)$ replaced by $V(n)$. So

$$(H\psi)(n) = \sum_{j:|j-n|=1} \psi(j) + V(n)\psi(n).$$

This $H$ is called a $d$-dimensional discrete Schrödinger operator.

We introduce two new spectra.

$$\sigma_{disc}(H) = \{E \in \sigma(H) \mid \text{ E is an isolated eigenvalue with finite multiplicity}\},$$

and

$$\sigma_{ess}(H) = \sigma(H) \backslash \sigma_{disc}(H).$$

The forth (and last) volume of "Methods of Modern Mathematical Physics" by Reed and Simon [37], called "Analysis of Operators", contains the following statement: "Spectral analysis of an operator $A$ concentrates on identifying the five sets $\sigma_{ess}(A)$, $\sigma_{disc}(A)$, $\sigma_{ac}(A)$, $\sigma_{sc}(A)$ and $\sigma_{pp}(A)$".

In this thesis, we will study the spectral properties of one-dimensional discrete limit-periodic Schödinger operators. We carry out our study within the framework of ergodic Schrödinger operators, a larger class to which our limit-periodic Schrödinger operators belong, so that the properties of ergodic Schrödinger operators automatically apply to limit-periodic ones. We will show that all of the three spectral types can appear within the class of limit-periodic Schrödinger operators. Besides, we will discuss the corresponding Lyapunov exponent, which plays an important role in describing time evolution of $\psi(t)$. In the region of pure point spectrum, we will exhibit the first almost periodic examples that are uniformly localized across the hull and the spectrum.

# Chapter 1

# An invitation to Limit-Periodic Schrödinger Operators

In the history of Schrödinger operator research, there are two exceptional landmark events. One is Anderson discovering the absence of diffusion for certain random lattice Hamiltonians in 1958, which finally led Anderson to win the Nobel prize in physics in 1977. Random Schrödinger operators are good models for such random lattices, like alloys, glasses, etc. Another event of importance is the discovery of quasicrystals by Shechtman in 1982, which caused a paradigm shift in crystallography. Shechtman won the Nobel prize in chemistry in 2011 for this surprising discovery. How do people mathematically study electronic properties of this disordered structure? Periodic Schrödinger operators, of which properties are well known, are used to describe ordered systems. While, for the study of those special disordered systems, quasi-periodic Schrödinger operators play an essential role. In this thesis, we will study limit-periodic Schrödinger operators, which are natural mathematical extensions of periodic ones, that is, the potentials can be uniformly approximated by periodic potentials. Also, quasi-periodic and limit-periodic Schrödinger operators both belong to the class of ergodic Schrödinger operators, where the potentials are defined dynamically, namely by sampling along the orbits of one or more ergodic transformations.

In this chapter, we will first introduce ergodic Schrödinger operators. Clearly, tools from dynamical systems are useful in the study of ergodic Schrödinger operators. In the subsequent sections we will mainly focus on the potential part, that is, a sequence in $\ell^\infty(\mathbb{Z})$, for it is the potential that makes Schrödinger operators really different. Most of the results in Section 1.3-1.5 have already been published in [24], where we exploit the connection between limit-periodic sequences and profinite groups.

## 1.1 Ergodic Schrödinger Operators

A discrete Schrödinger operator $H$ acting in $\ell^2(\mathbb{Z}^d)$ is given by

$$(H\psi)(i) = \sum_{j;|i-j|=1} \psi(j) + V(j)\psi(j). \tag{1.1}$$

We call $\{V(j)\}_{j\in\mathbb{Z}^d}$ a potential. In a random disordered system, the potential $\{V(j)\}_{j\in\mathbb{Z}^d}$ forms a random field, i.e., for any $j \in \mathbb{Z}^d$, $V(j)$ is a random variable on some probability space $(\Omega, \mathcal{F}, P)$. Without loss of generality, we can assume

$$\Omega = \prod_{\mathbb{Z}^d} S,$$

where $S$ is a Borel subset of $\mathbb{R}$ and $\mathcal{F}$ the $\sigma$-algebra generated by the cylinder sets of the form $\{\omega \mid \omega_{i_1} \in A_{i_1}, \omega_{i_2} \in A_{i_2}, \cdots, \omega_{i_n} \in A_{i_n}\}$ for $i_1, \cdots, i_n \in \mathbb{Z}^d$ and $A_{i_1}, \cdots, A_{i_n}$ are Borel sets in $S$. We define the shift operators $T^i$ on $\Omega$ by

$$(T^i(\omega))_j = \omega(j-i).$$

For a fixed $\omega$, the operator $H_\omega$ is nothing but a discrete Schrödinger operator with a certain potential. The point of random potentials is that we are no longer interested in properties of $H_\omega$ for a fixed $\omega$, but in properties suitable for a full measure set of $\omega$'s.

A probability measure $P$ on $\Omega$ is called *stationary* if $P(T^i(A)) = P(A)$ for any $A \in \mathcal{F}$ and $i \in \mathbb{Z}^d$. A stationary probability measure is called *ergodic* if any shift invariant set $A$, i.e., a set $A$ with $P(T^i(A)\Delta A) = 0$ for all $i \in \mathbb{Z}^d$, has probability $P(A) = 0$ or 1. Then $\{V(j)\}_{j\in\mathbb{Z}^d}$ can always be realized on the above probability space in such a way that $V_\omega(j) = \omega(j)$. $V$ is called stationary (ergodic), if the corresponding probability measure $P$ is stationary (ergodic). We introduce the celebrated ergodic theorem here.

**Theorem 1.1** (Ergodic Theorem, Birkhoff)**.** *Let $T$ be an ergodic transformation of a probability space $(\Omega, \mathcal{F}, P)$, and let $f : \Omega \to \mathbb{R}$ be a real measurable function. Then for almost every $\omega \in \Omega$ we have*

$$\lim_{n\to\infty} \frac{1}{n} \sum_{i=1}^{n} f(T^i(\omega)) = \int_\Omega f(\omega)dP.$$

Let's move to our target.

**Definition 1.2.** A one dimensional discrete ergodic Schrödinger operator acting in $\ell^2(\mathbb{Z})$ is given by

$$(H_{f,T}^{\omega}\psi)(i) = \psi(j+1) + \psi(j-1) + V_{\omega}(j)\psi(j) \tag{1.2}$$

with

$$V_{\omega}(j) = f(T^j(\omega)),$$

where $\omega$ belongs to a compact space $\Omega$, $T : \Omega \to \Omega$ a homeomorphism preserving an ergodic Borel probability measure $\mu$ and $f : \Omega \to \mathbb{R}$ a bounded measurable function.

Nice results about ergodic Schrödinger operators are the following.

**Theorem 1.3** (Pastur)**.** *For ergodic Schrödinger operators $H_{f,T}^{\omega}$, there exists a set $\Sigma \subset \mathbb{R}$ such that*

$$\sigma(H_{f,T}^{\omega}) = \Sigma \quad P-a.s.$$

*and*

$$\sigma_{dist}(H_{f,T}^{\omega}) = \emptyset \quad P-a.s.$$

**Theorem 1.4** (Kunz-Souillard)**.** *For ergodic Schrödinger operators $H_{f,T}^{\omega}$, there exists sets $\Sigma_{ac}, \Sigma_{sc}, \Sigma_{pp} \subset \mathbb{R}$ such that*

$$\sigma_{ac}(H_{f,T}^{\omega}) = \Sigma_{ac} \quad P-a.s.$$

$$\sigma_{sc}(H_{f,T}^{\omega}) = \Sigma_{sc} \quad P-a.s.$$

*and*

$$\sigma_{pp}(H_{f,T}^{\omega}) = \Sigma_{pp} \quad P-a.s.$$

One can find proofs of the above theorems in Simon et al's book [40]. So it is often beneficial to study the operators $\{H_{f,T}^{\omega}\}_{\omega\in\Omega}$ as a family, as opposed to a collection of individual operators.

In this thesis, we focus on one dimensional discrete Schrödinger operators. As a difference equation, what makes the one dimensional case accessible is that, for a fixed $E$, a solution of

$$(H_{f,T}^{\omega} - E)\psi = 0 \tag{1.3}$$

is totally determined by its value at two successive points. More precisely, $\psi : \mathbb{Z} \to \mathbb{R}$ is a solution of (1.3) if and only if

$$\begin{pmatrix} \psi(i+1) \\ \psi(i) \end{pmatrix} = \begin{pmatrix} E - f(T^i(\omega)) & -1 \\ 1 & 0 \end{pmatrix} \begin{pmatrix} \psi(i) \\ \psi(i-1) \end{pmatrix}. \tag{1.4}$$

Write $S_i = \begin{pmatrix} E - f(T^i(\omega)) & -1 \\ 1 & 0 \end{pmatrix}$ and $A_n^{(E,T,f)}(\omega) = S_n \dots S_1$. We call $A_n^{(E,T,f)}(\omega)$ a *transfer matrix*. The Lyapunov exponent is defined as

$$L(E,T,f) = \lim_{n\to\infty} \frac{1}{n} \int_\Omega \log \|A_n^{(E,T,f)}(\omega)\| \, d\mu(\omega), \tag{1.5}$$

which is to measure the growth rate of the solution of (1.3).

We would like to explain why the limit of (1.5) exists. A sequence $\{F_n\}_{n\in\mathbb{Z}}$ of random variables is called a *subadditive* process if

$$F_{n+m}(\omega) \le F_n(\omega) + F_m(T^n(\omega)),$$

where $T$ is a measure preserving transformation.

**Theorem 1.5** (Subadditive Ergodic Theorem, Kingman)**.** *If $F_n$ is a subadditive process satisfying $\mathbb{E}(|F_n|) < \infty$ for each $n$, and $\Gamma(F) := \inf \mathbb{E}(F_n/n) > -\infty$, then $F_n(\omega)/n$ converges almost surely. If, furthermore, $T$ is ergodic, then*

$$\lim_{n\to\infty} F_n(\omega)/n = \Gamma(F)$$

*almost surely.*

Let $F_n(\omega) = \log \|A_n^{(E,T,f)}(\omega)\|$. It is easy to see that

$$F_{n+m}(\omega) \le F_n(\omega) + F_m(T^n(\omega)),$$

so that $F_n(\omega)$ is a subadditive process. We also have

$$\begin{aligned} \frac{1}{n}\mathbb{E}(|F_n|) &= \frac{1}{n}\mathbb{E}(\log \| \prod_n^1 S_i \|) \\ &\le \frac{1}{n}\mathbb{E}(\sum_i^n \log \| S_i \|) \\ &= \mathbb{E}(\log \| S_0 \|), \end{aligned}$$

where we use the stationarity of $V_\omega(i)$. Since $V_\omega(i) = f(T^i(\omega))$ is bounded, we have $\mathbb{E}(\log \| S_0 \|) < \infty$. In addition, we have $\mathbb{E}(F_n/n) \geq 0 > -\infty$. By the subadditive ergodic theorem, we have

$$\lim_{n\to\infty} \frac{1}{n} \log \|A_n^{(E,T,f)}(\omega)\| = \inf_{n>0} \frac{1}{n}\mathbb{E}(\log \|A_n^{(E,T,f)}(\omega)\|) \quad a.s.$$

By the dominated convergence theorem, $\lim_{n\to\infty} \frac{1}{n} \int_\Omega \log \|A_n^{(E,T,f)}(\omega)\| \, d\mu(\omega)$ exists.

We introduce three useful properties about the Lyapunov exponent to close this section.

**Theorem 1.6** (Craig-Simon)**.** *$L(E,T,F)$ is a subharmonic function of $E \in \mathbb{R}$.*

**Theorem 1.7** (lshii-Pastur-Kotani)**.** *Suppose that $V_\omega$ is a bounded ergodic process. Then*

$$\sigma_{a.c.}(H^\omega_{f,T}) = \overline{\{E \mid L(E,f,T) = 0\}}^{ess}.$$

*Here the essential closure is defined by*

$$\bar{A}^{ess} = \{a \mid |A \cap (a-\varepsilon, a+\varepsilon)| > 0 \text{ for all } \varepsilon > 0\}.$$

**Theorem 1.8** (Kotani)**.** *If $V_\omega(n)$ is nondeterministic, then $L(E,f,T) > 0$ for Lebesgue-almost all $E \in \mathbb{R}$. Thus, $\sigma_{ac}(H^\omega_{f,T}) = \emptyset$.*

## 1.2 Almost Periodic Sequences

In this section we want to review some basic properties of almost periodic sequences [1], for limit-periodic and quasi-periodic sequences are both almost periodic.

**Definition 1.9.** Let $T : \ell^\infty(\mathbb{Z}) \to \ell^\infty(\mathbb{Z})$ be the shift operator, $(T(V))(n) = V(n+1)$. A two-sided sequence $V \in \ell^\infty(\mathbb{Z})$ is called almost periodic if the closure of its $T$-orbit is compact. The compact space obtained by taking the closure of its $T$-orbit is called the hull of $V$, denoted by $\Omega_V$.

[1] In this chapter we are using "sequence" and in the second chapter we will use "potential" instead.

The compact space $\Omega_V$ has the natural abelian group structure. We define the group operation as $T^i(V) \cdot T^j(V) = T^{i+j}(V)$. We let $W$ and $\tilde{W}$ be the limits of $T^{i_n}(V)$ and $T^{j_n}(V)$ respectively. A simple calculation shows

$$\| T^{i_n+j_n}(V) - T^{i_m+j_m}(V) \| \leq \| T^{i_n}(V) - T^{i_m}(V) \| + \| T^{j_n}(V) - T^{j_m}(V) \| .$$

This implies that the limit of $T^{i_n+j_n}(V)$ exists, and denoted by $W \cdot \tilde{W}$. So the group operation is well defined. Also, $\mathbb{Z}$ can be embedded to be a dense supgroup of $\Omega_V$ by $n \to T^n(V)$.

$\Omega_V$ as a compact abelian topological group has a natural normalized Haar measure $d\mu_V$, for which the shift operator $T$ is an ergodic transformation. More precisely, for any $f \in C(\Omega_V, \mathbb{R})$, we have

$$\int_{\Omega_V} f(\omega) d\mu_V = \lim_{n\to\infty} \frac{1}{2n} \sum_{-n}^{n} f(T^n(W)), a.s.\ W \in \Omega_V,$$

which follows from the ergodic theorem. When $T$ is unique ergodic, the above equality is successful for every $W \in \Omega_V$.

We then see that every almost periodic sequence has an average and that this average is Haar measure on the hull of that sequence. Conversely, any continuous function $f$ on $\Omega_V$ defines an almost periodic sequence $\tilde{V}$ by $\tilde{V}(i) = f(T^i(V))$ since $\{\tilde{V}(i)\}_{i\in\mathbb{Z}}$ lies in the image of $\Omega_V$ under a continuous map; the hull of $\tilde{V}$ is naturally a quotient group of $\Omega_V$.

The dual group, $\hat{\Omega}_V$, of characters on $\Omega_V$ is naturally a topological subgroup of the circle group $\mathbb{T}$, the multiplicative group of all complex numbers of absolute value 1, since $\mathbb{Z}$ can be embedded into $\Omega_V$ as the dense subgroup and the dual group of $\mathbb{Z}$ is $\mathbb{T}$. By taking inverse image under the map $\mathbb{R} \to \mathbb{T}, \alpha \to e^{i\alpha}$, we obtain a subgroup $F_V$ of $\mathbb{R}$, called the *frequency module* of $V$. $F_V$ is countable since $\Omega_V$ has a countable dense subset orb$(V)$, and it is a $\mathbb{Z}$-module. By [42, Theorem 5.13A.1.], we have that $V$ is almost periodic if and only if $V$ is a uniform limit of finite sums of the form

$$Q_N(n) = \sum_{j=1}^{N} c_j e^{2\pi i \alpha_j^{(N)} n}$$

for $\alpha_1^{(N)}, \cdots, \alpha_N^{(N)} \in \mathbb{R}/\mathbb{Z}$. We also have the following theorem.

**Theorem 1.10.** [3, Theorem A.1.1] *The frequency module $F_V$, is the $\mathbb{Z}$-module generated by*

$$\{\alpha : \lim_{n\to\infty} \frac{1}{2n} \sum_{k=-n}^{n} V(k)e^{-ik\alpha} \neq 0, \alpha \in \mathbb{R}\}.$$

Thus, every countable module in $\mathbb{Z}$ is the frequency module of some almost periodic sequence.

**Definition 1.11.** $V \in \ell^\infty(\mathbb{Z})$ is called periodic if there exists $p \in \mathbb{Z}^+$ so that $V(n) = V(n+p)$ for every $n \in \mathbb{Z}$.

The hull of a periodic sequence is a finite set. There are two other classes of almost periodic sequences of special interest.

**Definition 1.12.** We say that $V$ is limit-periodic if and only if $V$ is a uniform limit of periodic sequences. We say that $V$ is quasi-periodic if and only if there exist numbers $\alpha_1, \alpha_2, \cdots, \alpha_m$ and a continuous function $f$ on $\mathbb{T}^m$ (the $m$-torus) with

$$V(n) = f([n\alpha_1], [n\alpha_2], \cdots, [n\alpha_m]),$$

with $[n\alpha_i]$ the fractional part of $n\alpha_i$ (thinking of $\mathbb{T}$ as $\mathbb{R}/\mathbb{Z}$).

If the $\alpha_i$ are rationally independent, then $\Omega_V$ is the image of $\mathbb{T}^m$ under $f$. The quasi-periodic one is a good model to study quasicystals. The hull of a limit-periodic sequence is a totally disconnected compact abelian topological group (see Proposition 1.20). So they are both almost periodic. One can tell these special classes apart by looking at the frequency module.

**Theorem 1.13.** [3, Theorem A.1.2] *Let $V$ be almost periodic. $V$ is quasi-periodic if and only if $F_V$ is finitely generated.*

**Theorem 1.14.** [3, Theorem A.1.3] *Let $V$ be almost periodic. $V$ is limit-periodic if and only if $F_V$ has the property that any $\alpha, \beta \in F_V$ have a common divisor in $F_V$.*

**Corollary 1.15.** [3, Corollary A.1.4] *If $V$ is both limit-periodic and quasi-periodic, then $V$ is periodic.*

**Corollary 1.16.** [3, Corollary A.1.5] *If $V$ is limit-periodic, then there exists a positive integer set $I_V = \{n_j\}$ satisfying $n_j | n_{j+1}$ such that*

$$V(k) = \sum_{j=1}^{\infty} P_j(k),$$

*where the $P_j \in \ell^\infty(\mathbb{Z})$ are $n_j$-periodic, and the convergence is uniform.*

We would like to state one more property to close this section.

**Theorem 1.17.** [3, Theorem A.2.2.] *Let $V$ be almost periodic. Then the spectrum of $-\Delta + V$ is purely essential, i.e. there are no isolated eigenvalues of finite multiplicity.*

## 1.3 Limit-Periodic Sequences

In [1], Avila gave a way to treat limit-periodic sequences by regarding limit-periodic sequences as generated by continuous functions along the orbits of a minimal translation of a Cantor group, allowing one to separate the base dynamics and the sampling function. Though this is quite standard in the quasi-periodic case, it is still new in the limit-periodic case.

**Definition 1.18.** A group $\Omega$ is called Cantor if it is a totally disconnected compact abelian topological group without isolated points.

**Definition 1.19.** Given a topological group $\Omega$, a map $T : \Omega \to \Omega$ is called a translation if $T(\omega) = \omega \cdot \omega_0$ for some $\omega_0 \in \Omega$, and $T$ is called minimal if the orbit $\{T^n(\omega) : n \in \mathbb{Z}\}$ of every $\omega \in \Omega$ is dense in $\Omega$.

For $V \in \ell^\infty(\mathbb{Z})$, let $\mathrm{orb}(V) = \{T^k(V) : k \in \mathbb{Z}\}$, and $\mathrm{hull}(V)$ be the closure of $\mathrm{orb}(V)$ in $\ell^\infty$-norm.

**Proposition 1.20.** [1, Lemma 2.1] *If $V$ is a limit-periodic sequence, then* $\mathrm{hull}(V)$ *is compact and has a unique topological group structure with identity $T^0(V) = V$ such that*

$$\phi : \mathbb{Z} \longrightarrow \mathrm{hull}(V), \quad k \longrightarrow T^k(V)$$

*is a homomorphism. Also, the group structure is abelian and there exist arbitrarily small compact open neighborhoods of $V$ in* $\mathrm{hull}(V)$ *that are finite index subgroups.*

Proposition 1.20 tells us that a limit-periodic sequence is almost periodic and furthermore that hull($V$) is totally disconnected. So hull($V$) is just a Cantor group that admits a minimal translation. We denote it by $\Omega_V$.

**Remark 1.21.** *Not every Cantor group has minimal translations. For example,*

$$\Omega = \prod_{j=0}^{\infty} \mathbb{Z}_2,$$

*where $\mathbb{Z}_2$ is a cyclic 2-group, is a Cantor group with the product topology, but it has no minimal translations, since $\omega + \omega$ is the identity for any $\omega \in \Omega$.*

In the last section, we see that an almost periodic sequence can be determined by its frequency module. Given a limit-periodic potential $V \in \ell^\infty(\mathbb{Z})$ which is not periodic, its frequency module $F_V$ is countable and infinitely generated (if it is finitely generated, $V$ is quasi-periodic or periodic, see Theorem 1.13). Denote the set of generators of $F_V$ by $G_V = \{2\pi\alpha_j\}$. $2\pi\alpha_1$ and $2\pi\alpha_2$ have a common divisor by Theorem 1.14. By Theorem 1.10,

$$a = \lim_{n\to\infty} \frac{1}{2n} \sum_{k=-n}^{n} V(k) e^{-i2k\pi\alpha_1}$$

and

$$b = \lim_{n\to\infty} \frac{1}{2n} \sum_{k=-n}^{n} V(k) e^{-i2k\pi\alpha_2}$$

are both non-zero. Choose a periodic potential $P \in \ell^\infty(\mathbb{Z})$ with

$$\| P - V \|_\infty \leq \frac{1}{2} \min(|a|, |b|).$$

It follows that

$$\lim_{n\to\infty} \frac{1}{2n} \sum_{k=-n}^{n} P(k) e^{-i2k\pi\gamma} \neq 0$$

for $\gamma = \alpha_1, \alpha_2$. So $\frac{2\pi}{h}$ divides both $2\pi\alpha_1$ and $2\pi\alpha_2$ where $h$ is the period [2] of $P$ (note that the frequency module of an $h$-periodic sequence is generated by $\frac{2\pi}{h}$). So there exists $n_1 \in \mathbb{Z}^+$ such that the greatest common divisor of $2\pi\alpha_1$ and $2\pi\alpha_2$ is $\frac{2\pi}{n_1}$, and the $\mathbb{Z}$-module generated by $\{2\pi\alpha_1, 2\pi\alpha_2\}$ is the $\mathbb{Z}$-module generated by $\{\frac{2\pi}{n_1}\}$. Similarly, for $\alpha_1, \alpha_2$ and $\alpha_3$ there exists a positive integer $n_2$ such that the $\mathbb{Z}$-module generated by $\{2\pi\alpha_1, 2\pi\alpha_2, 2\pi\alpha_3\}$ is the $\mathbb{Z}$-module generated by $\{\frac{2\pi}{n_2}\}$. Clearly, $n_1|n_2$. By induction, we can find an infinite positive integer set

$$I_V = \{n_j\} \subset \mathbb{Z}^+$$

such that $n_j|n_{j+1}$ and $F_V$ is generated by $G_V = \left\{\frac{2\pi}{n_j} : n_j \in I_V\right\}$.

[2] In this thesis, the period of a periodic sequence is considered as the minimal period.

We call $I_V$ a *frequency integer set* of $V$. In general, we say that $I_V$ is a frequency integer set of $V$ if one can obtain the frequency module $F_V$ from $I_V$. When $V$ is limit-periodic, for any $n, m \in I_V$, one always has $n|m$ or $m|n$. If $V(k) = \sum_{j=1}^{\infty} P_j(k)$, where the $P_j$ are $n_j$-periodic, then it is easy to conclude that $I_V = \{n_j\}$ is a frequency integer set of $V$. Clearly, $V$ may have other frequency integer sets. By the following theorem, we shall see that there exists a unique maximal frequency integer set $M_V$, in the sense that every frequency integer set $I_V$ is contained in $M_V$.

**Theorem 1.22.** *Given a limit-periodic sequence $V \in \ell^\infty(\mathbb{Z})$ with an infinite frequency integer set $I_V$, for any limit-periodic sequence $\tilde{V} \in \ell^\infty(\mathbb{Z})$ with a frequency integer set $I_{\tilde{V}}$, if $I_{\tilde{V}}$ is an infinite subset of $I_V$, then we have $\Omega_V \cong \Omega_{\tilde{V}}$.*

*Proof.* Write $I_{\tilde{V}} \equiv \{n_{j_k}\} \subset I_V \equiv \{n_j\}$. Define a homomorphism by $\varphi : F_{\tilde{V}} \to F_V$, $\varphi(\frac{2\pi}{n_{j_k}}) = \frac{2\pi}{n_{j_k}}$. Clearly, $\varphi$ is injective since $G_{\tilde{V}} \subset G_V$. For any $\frac{2\pi}{n_t} \notin G_{\tilde{V}}$, choose $n_{j_k} > n_t$, and then there is $t_k$ such that $n_{j_k} = t_k n_t$ (by the property of $I_V$). We have $\varphi(t_k \frac{2\pi}{n_{j_k}}) = t_k \frac{2\pi}{n_{j_k}} = \frac{2\pi}{n_t}$, implying that $\varphi$ is also surjective. Thus, $F_V \cong F_{\tilde{V}}$. Clearly, $\hat{\Omega}_V \cong \hat{\Omega}_{\tilde{V}}$, since $\hat{\Omega}_V = \{e^{i\alpha} : \alpha \in F_V\}$ and $\hat{\Omega}_{\tilde{V}} = \{e^{i\beta} : \beta \in F_{\tilde{V}}\}$. By the Pontryagin duality theorem, $\Omega_V \cong \Omega_{\tilde{V}}$. □

Theorem 1.22 gives us a way to find the maximal frequency integer set of $V$. If $I_V$ is a frequency integer set of $V$, one can add positive integers into $I_V$ to get $M_V = \{m_j\} \subset \mathbb{Z}^+$ such that $m_{j+1}/m_j$ is a prime number. Clearly, $I_V \subset M_V$ and $M_V$ is still a frequency integer set of $V$. The uniqueness of $M_V$ follows from Theorem 1.22.

**Theorem 1.23.** *Given limit-periodic potentials $V, \tilde{V} \in \ell^\infty(\mathbb{Z})$ with frequency integer sets $I_V$ and $I_{\tilde{V}}$ respectively, $\Omega_V \cong \Omega_{\tilde{V}}$ if and only if for any $n_i \in I_V$ there exists $m_j \in I_{\tilde{V}}$ such that $n_i|m_j$, and vice versa.*

*Proof.* If $S_V, S_{\tilde{V}}$ are both finite sets, $V, \tilde{V}$ are periodic, and then the statement follows readily. We assume that $I_V, I_{\tilde{V}}$ are infinite sets. First assume that for any $n_i \in I_V$ there exists $m_j \in I_{\tilde{V}}$ such that $n_i | m_j$ and vice versa. Then, for $n_{i_1} \in I_V$, there exists $m_{j_1} \in I_{\tilde{V}}$ such that $n_{i_1} | m_{j_1}$. Consider $m_{j_1} \in I_{\tilde{V}}$, and similarly, there exists $n_{i_2} \in I_V$ with $m_{j_1} | n_{i_2}$. By induction, we get infinite subsets $L_V \equiv \{n_{i_k}\} \subset I_V$ and $L_{\tilde{V}} \equiv \{m_{j_k}\} \subset I_{\tilde{V}}$ such that $H \equiv L_V \cup L_{\tilde{V}}$ can be a frequency integer set of some limit-periodic potential $V'$ (see Corollary 1.16). By Theorem 1.22, we conclude that $\Omega_V \cong \Omega_{V'} \cong \Omega_{\tilde{V}}$.

Conversely, assume that such a condition is not successful. Without loss of generality, suppose that for a fixed $n_1 \in I_V$, $n_1 \nmid m_j$ for any $m_j \in I_{\tilde{V}}$. If $\Omega_V \cong \Omega_{\tilde{V}}$, by the Pontryagin duality theorem, we have $\hat{\Omega}_V \cong \hat{\Omega}_{\tilde{V}}$. So there is an isomorphism $\phi : \hat{\Omega}_V \to \hat{\Omega}_{\tilde{V}}$ with $\phi(e^{2\pi i}) = e^{2\pi i}$ (the identity must be mapped to the identity). Write

$$\phi(e^{\frac{2\pi}{n_1}i}) = e^{\sum_{j=1}^{t_1} k_j \frac{2\pi i}{m_j}} = e^{q_1 \frac{2\pi i}{m_{t_1}}},$$

and we get

$$\phi(e^{m_{t_1}\frac{2\pi}{n_1}i}) = e^{q_1 m_{t_1} \frac{2\pi i}{m_{t_1}}} = e^{2q_1\pi i} = 1,$$

which is impossible since $m_{t_1}/n_1$ is not an integer. Thus, $\Omega_V \ncong \Omega_{\tilde{V}}$. The proof is complete. □

By Theorem 1.23, we conclude the following classification theorem about limit-periodic sequences.

**Theorem 1.24.** *Given two limit-periodic potentials $V, \tilde{V} \in \ell^\infty(\mathbb{Z})$, we have $\Omega_V \cong \Omega_{\tilde{V}}$ if and only if $V$ and $\tilde{V}$ have the same maximal frequency integer set.*

If $V \in \ell^\infty(\mathbb{Z})$ is periodic, hull$(V)$ is a finite cyclic group. If $V$ is limit-periodic but not periodic, the hull of a limit-periodic sequence is a Cantor group that admits a minimal translation. We shall see that hull$(V)$ is just a procyclic group.

## 1.4 Profinite Groups

In this section we will discuss profinite groups. We will see that Cantor groups actually belong to the class of profinite groups.

A topological profinite group is by definition an *inverse limit* of finite topological groups. Now let's introduce the related definitions. A *directed* set is a partially ordered set $I$ such that for all $i_1, i_2 \in I$ there is an element $j \in I$ for which $i_1 \leq j$ and $i_2 \leq j$.

**Definition 1.25.** An inverse system $(X_i, \phi_{ij})$ of topological groups indexed by a directed set $I$ consists of a family $(X_i \mid i \in I)$ of topological groups and a family $(\phi_{ij} : X_j \to X_i \mid i, j \in I, i \leq j)$ of continuous homomorphisms such that $\phi_{ii}$ is the identity map $\mathrm{id}_{X_i}$ for each $i$ and $\phi_{ij}\phi_{jk} = \phi_{ik}$ whenever $i \leq j \leq k$.

For example, let $I = \mathbb{Z}^+$ with the usual order, let $p$ be a prime, and let $X_i = \mathbb{Z}/p^i\mathbb{Z}$ for each $i$, and for $j \geq i$ let $\phi_{ij} : X_j \to X_i$ be the map defined by

$$\phi_{ij}(n + p^j\mathbb{Z}) = n + p^i\mathbb{Z}$$

for each $n \in \mathbb{Z}$. Then $(X_i, \phi_{ij})$ is an inverse system of finite topological groups.

Let $(X_i, \phi_{ij})$ be an inverse system of topological groups and let $Y$ be a topological group. We call a family $(\phi_i : Y \to X_i \mid i \in I)$ of continuous homomorphisms *compatible* if $\phi_{ij}\phi_j = \phi_i$ whenever $i \leq j$; that is, the following diagram

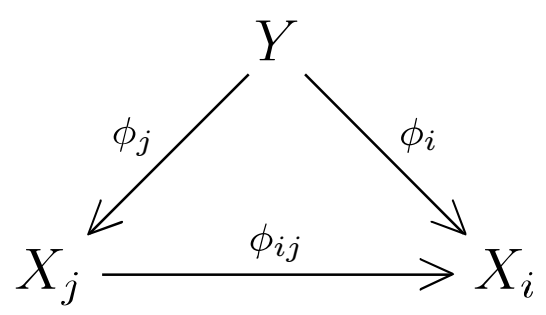

is commutative.

**Definition 1.26.** An inverse limit $(X, \phi_i)$ of an inverse system $(X_i, \phi_{ij})$ of topological groups is a topological group together with a compatible family $(\phi_i : X \to X_i)$ of continuous homomorphisms with the following universal property: whenever $(\varphi_i : Y \to X_i)$ is a compatible family of continuous homomorphisms from a topological group $Y$, there is a unique continuous homomorphism $\varphi : Y \to X$ such that $\phi_i\varphi = \varphi_i$ for each $i$.

That is, there is a unique continuous homomorphism $\varphi$ such that the following diagram

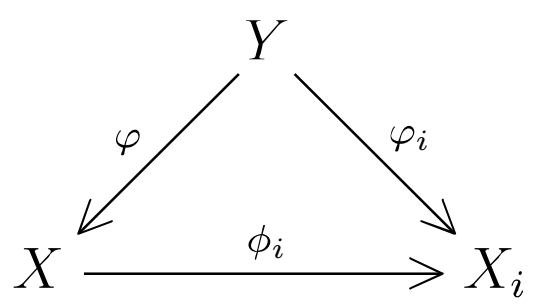

is commutative.

**Proposition 1.27.** [46, Proposition 1.1.4] *Let $(X_i, \phi_{ij})$ be an inverse system of topological groups, indexed by $I$.*

(1). *There exists an inverse limit $(X, \phi_i)$ of $(X_i, \phi_{ij})$, for which $X$ is a topological group and the maps $\phi_i$ are continuous homomorphisms.*

(2). *If $(X^{(1)}, \phi_i^{(1)})$ and $(X^{(2)}, \phi_i^{(2)})$ are inverse limits of the inverse system, then there is an isomorphism $\bar{\phi} : X^{(1)} \to X^{(2)}$ such that $\phi_i^{(2)}\bar{\phi} = \phi_i^{(1)}$ for each $i$.*

(3). *Write $G = \prod_{i \in I} X_i$ with the product topology and for each $i$ write $\pi_i$ for the projection map from $G$ to $X_i$. Define*

$$X = \{c \in G : \phi_{ij}\pi_j(c) = \pi_i(c) \ for\ all\ i, j\ with\ j \geq i\}$$

*and $\phi_i = \pi_i|_X$ for each $i$. Then $(X, \phi_i)$ is an inverse limit of $(X_i, \phi_{ij})$.*

The result above shows that the inverse limit of an inverse system $(X_i, \phi_{ij})$ exists and is unique up to isomorphism. We also have the following important characterization of profinite groups.

**Proposition 1.28.** [46, Corollary 1.2.4] *Let $X$ be a topological group. The following are equivalent:*

(1). *$X$ is profinite, i.e., it is an inverse limit of an inverse system;*

(2). *$X$ is isomorphic to a closed subgroup of a product group of finite groups;*

(3). *$X$ is compact and $\bigcap(N \mid N \lhd_O X) = 1$ ($\lhd_O$ means that $N$ is open and normal);*

(4). *$X$ is compact and totally disconnected.*

By the above proposition, we see that *a Cantor group is an abelian profinite group without isolated points.*

**Proposition 1.29.** *Assume that the directed set is $I = \mathbb{Z}^+$ with the usual order. For an inverse system $(X_i, \phi_{ij})_{j \geq i}$ with the inverse limit $(X, \phi_i)$, every non-finite sub-inverse system still has the same inverse limit $(X, \phi_{i_k})$ up to isomorphism.*

*Proof.* Consider a non-finite sub-inverse system $(X_{i_k}, \phi_{i_k i_t})_{t \geq k}$. Assume that $(X', \phi^{(1)}_{i_k})$ is the inverse limit. Obviously, $(X, \phi_{i_k})$ is compatible with $(X_{i_k}, \phi_{i_k i_t})_{t \geq k}$, so there is a unique homomorphism $\phi^{(1)} : X \to X'$ such that $\phi_{i_k} = \phi^{(1)}_{i_k}\phi^{(1)}$.

For any $X_q$ not in the sub-inverse system, choose $i_k > q$. We have that $\phi_{q i_k} : X_{i_k} \to X_q$ and $\phi^{(1)}_q = \phi_{q i_k}\phi^{(1)}_{i_k} : X' \to X_q$ are homomorphisms. We will prove that $(X', \phi^{(1)}_i)$ is compatible with $(X_i, \phi_{ij})_{j \geq i}$, for which it suffices to show that the following diagram:

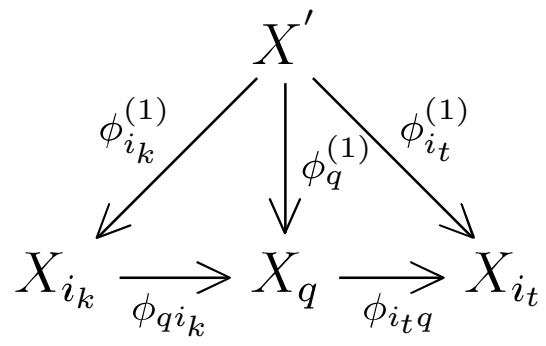

is commutative. The left half of the above diagram follows from the definition of $\phi^{(1)}_q$. The right half follows from $\phi^{(1)}_{i_t} = \phi_{i_t q}\phi_{q i_k}\phi^{(1)}_{i_k} = \phi_{i_t q}\phi^{(1)}_q$. So $(X', \phi^{(1)}_i)$ is compatible with $(X_i, \phi_{ij})_{j \geq i}$ and there is a unique homomorphism $\phi^{(2)} : X' \to X$ such that $\phi^{(1)}_i = \phi_i \phi^{(2)}$.

By the universal property for $(X', \phi^{(1)}_{i_k})$, there is only one map $F : X' \to X'$ with the property that $\phi^{(1)}_{i_k} F = \phi^{(1)}_{i_k}$ for each $i$. However, both $\phi^{(1)}\phi^{(2)}$ and $\mathrm{id}_{X'}$ have this property, so $\phi^{(1)}\phi^{(2)} = \mathrm{id}_{X'}$. Similarly, for $\phi^{(2)}\phi^{(1)} : X \to X$ we have $\phi_{i_k}\phi^{(2)}\phi^{(1)} = \phi_{i_k}$. For any $\phi_j$, $\phi_j = \phi_{j i_k}\phi_{i_k}$ when $i_k > j$, so we have $\phi_j \phi^{(2)}\phi^{(1)} = \phi_{j i_k}\phi_{i_k}\phi^{(2)}\phi^{(1)} = \phi_{j i_k}\phi_{i_k} = \phi_j$, which implies $\phi^{(2)}\phi^{(1)} = \mathrm{id}_X$. Thus, $\phi^{(1)} : X \to X'$ is an isomorphism. □

**Proposition 1.30.** [46, Proposition 1.1.7] *Let $G$ be a compact Hausdorff totally disconnected space. Then $G$ is the inverse limit of its discrete quotient spaces.*

We interpret a *class* in the usual sense that it is closed with respect to taking isomorphic images. Let $\mathcal{C}$ be some class of finite groups. We call a group $F$ a $\mathcal{C}$-group if $F \in \mathcal{C}$, and $G$ is called a pro-$\mathcal{C}$ group if it is an inverse limit of $\mathcal{C}$-groups. We say that $\mathcal{C}$ is closed for quotients (resp. subgroups) if every quotient group (resp. subgroup) of a $\mathcal{C}$-group is also a $\mathcal{C}$-group. Similarly, we say that $\mathcal{C}$ is closed for direct products if $F_1 \times F_2 \in \mathcal{C}$ whenever $F_1 \in \mathcal{C}$ and $F_2 \in \mathcal{C}$. For example, for the class of finite $p$-groups where $p$ is a fixed prime, an inverse limit of finite $p$-groups is called a *pro-p* group; for the class of finite cyclic groups, an inverse limit of finite cyclic groups is called a *procyclic* group.

The next result describes how a given profinite group, its subgroups and quotient groups, can be represented explicitly as inverse limits.

**Proposition 1.31.** [46, Theorem 1.2.5] (1). *Let $G$ be a profinite group. If $I$ is a filter base of closed normal subgroups of $G$ such that $\bigcap(N \mid N \in I) = 1$, then*

$$G \cong \varprojlim_{N \in I} G/N.$$

*Moreover*

$$H \cong \varprojlim_{N \in I} H/(H \cap N)$$

*for each closed subgroup $H$, and*

$$G/K \cong \varprojlim_{N \in I} G/KN$$

*for each closed normal subgroup $K$.*
(2). *If $\mathcal{C}$ is a class of finite groups which is closed for subgroups and direct products, then closed subgroups, direct products and inverse limits of pro-$\mathcal{C}$ groups are pro-$\mathcal{C}$ groups. If in addition $\mathcal{C}$ is closed for quotients, then quotient groups of pro-$\mathcal{C}$ groups by closed normal subgroups are pro-$\mathcal{C}$ groups.*

## 1.5 Why Profinite Groups

Let's discuss the connection between limit-periodic sequences and profinite groups. Please refer to [24] for more details.

### 1.5.1 Generate Limit-Periodic Sequences.

First, let's see how to generate limit-periodic sequences from a Cantor group which admits a minimal translation. (Note that such a group is just a profinite group. Moreover, we shall see that such a group is procyclic).

**Lemma 1.32.** [1, Lemma 2.2.] *Given a Cantor group $\Omega$ and a minimal translation $T$, for each $f \in C(\Omega, \mathbb{R})$, define $F : \Omega \to \ell^\infty(\mathbb{Z})$, $F(\omega) = (f(T^n(\omega)))_{n\in\mathbb{Z}}$. We have that $F(\omega)$ is limit-periodic and $F(\Omega) = \mathrm{hull}(F(\omega))$ for every $\omega \in \Omega$.*

By Proposition 1.20, we know that $\mathrm{hull}(F(e))$ ($e$ is the identity of $\Omega$) is a finite cyclic group or it is a Cantor group with the unique group structure: $T^0(F(e)) = F(e)$ is the identity element, and $T^i(F(e)) \cdot T^j(F(e)) = T^{i+j}(F(e))$.

Since $C(\Omega, \mathbb{R})$ can generate a class of limit-periodic sequences, we will get a certain class of topological groups by taking the hulls of these limit-periodic sequences. The group in this class is a Cantor group or a finite cyclic group. What is the relation between this class of topological groups and the original Cantor group $\Omega$? (Note that we will adopt the notation $F(e)$ throughout this section, that is, $F(e) = (f(T^n(e)))_{n\in\mathbb{Z}}$ as in Lemma 1.32.)

**Lemma 1.33.** *There exists an $f \in C(\Omega, \mathbb{R})$ such that* $\mathrm{hull}(F(e)) \cong \Omega$.

*Proof.* It is easy to see that for any $f \in C(\Omega, \mathbb{R})$, $\mathrm{hull}(F(e))$ is a quotient group of $\Omega$. So it suffices to prove that there exists some $f \in C(\Omega, \mathbb{R})$ such that $\ker(\phi) = \{e\}$. Clearly, a Cantor group is metrizable (recall that any separable compact space is metrizable). Introduce a metric on $\Omega$ compatible with the topology. Define a function $f : \Omega \to \mathbb{R}$ by $f(\omega) = \mathrm{dist}(e, \omega)$. Clearly, $f$ is continuous, so there is a corresponding $F$ (defined as in Lemma 1.32) $: \Omega \to \ell^\infty(\mathbb{Z})$ such that $\mathrm{hull}(F(e))$ is a quotient group of $\Omega$. Consider $\phi : \Omega \to \mathrm{hull}(F(e)), \phi(\omega) = F(\omega)$. If $F(\omega) = F(e)$, then $f(\omega) = f(e)$, that is, $\mathrm{dist}(e, \omega) = \mathrm{dist}(e, e) = 0$, implying $\omega = e$ and $\ker(\phi) = \{e\}$. □

**Theorem 1.34.** *Given a Cantor group $\Omega$ and a minimal translation $T$, for each $f \in C(\Omega, \mathbb{R})$,* $\mathrm{hull}(F(e))$ *is a Cantor group or a finite cyclic group, and $C(\Omega, \mathbb{R})$ can generate a class of topological groups. Then there is a one-to-one correspondence between this class of topological groups and quotient groups of $\Omega$.*

*Proof.* We know that groups in this class are quotient groups of $\Omega$. Let's show the converse direction. By the definition, a quotient group of a Cantor group is still Cantor or finitely cyclic. Given a quotient group $\Omega_0$ and a quotient homomorphism $q : \Omega \to \Omega_0$, we claim that $T$ will induce a minimal translation $T_0$ on $\Omega_0$ such that $T_0([\omega]) = q(T(\omega)), [\omega] \in \Omega_0$. For writing convenience, we assume that the group operation is addition and $T(\omega) = \omega + \omega_1$. If $[\omega] = [\omega']$, then $T_0([\omega]) = q(T(\omega)) = q(\omega_1 + \omega) = q(\omega_1) + q(\omega) = [\omega_1] + [\omega] = [\omega_1] + [\omega'] = T_0([\omega'])$, which gives that $T_0$ is a translation by $[\omega_1]$. That $T_0$ is minimal follows from the fact that $T$ is minimal and $q$ is continuous. By Lemma 1.33, we know that for $\Omega_0$ and $T_0$ there exists some $f_0 \in C(\Omega_0, \mathbb{R})$ such that $\text{hull}((f_0(T_0^n([e])))_{n\in\mathbb{Z}}) \cong \Omega_0$. Let $f = f_0 \circ q$. Clearly, $f \in C(\Omega, \mathbb{R})$ and the following diagram is commutative. So we have $\text{hull}(F(e)) \cong \text{hull}((f_0(T_0^n([e])))_{n\in\mathbb{Z}}) \cong \Omega_0$. The proof is complete. □

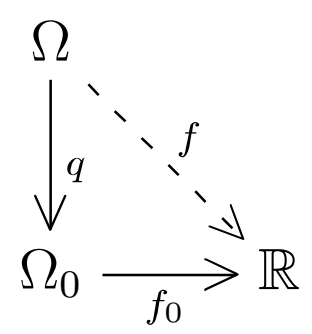

We also have

**Theorem 1.35.** *Given a Cantor group* $\Omega$ *and a minimal translation* $T$*, for any limit-periodic sequence* $V \in \ell^\infty(\mathbb{Z})$ *satisfying* $\text{hull}(V) \cong \Omega$*, there is an* $f \in C(\Omega, \mathbb{R})$ *such that* $f(T^n(e)) = V_n$ *for every* $n \in \mathbb{Z}$*.*

*Proof.* By Lemma 1.33 we have $\tilde{f} \in C(\Omega, R)$ such that $\text{hull}(\tilde{F}(e)) \cong \Omega$ (note that $\tilde{F}(e) = (\tilde{f}(T^n(e)))_{n\in\mathbb{Z}}$). Since $\text{hull}(V) \cong \Omega$, we have a continuous isomorphism $h : \text{hull}(\tilde{F}(e)) \to \text{hull}(V)$ with $h(\tilde{F}(e)) = V$.

Clearly, for $T^{n_k}(e) \in \Omega$ we have $h(\tilde{F}(T^{n_k}(e))) = \tilde{T}^{n_k}(V)$ ($\tilde{T}$ is the left shift operator of $V$ ) since $\tilde{F}(T^{n_k}(e)) = \tilde{T}^{n_k}(\tilde{F}(e))$. If $\lim_{k\to\infty} T^{n_k}(e) = \omega$, then $h(\tilde{F}(\omega)) = \lim_{k\to\infty} \tilde{T}^{n_k}(V)$, where the limit exists since $h$ and $\tilde{F}$ are both continuous. Define $f$ by $f(T^n(e)) = \tilde{T}^n(V)_0 = V_n$. We extend $f$ to the whole $\Omega$ by $f(\omega) = \lim_{k\to\infty} \tilde{T}^{n_k}(V)_0$ if $\omega = \lim_{k\to\infty} T^{n_k}(e)$. By the previous analysis, $f$ is well defined and continuous. So there is an $f \in C(\Omega, \mathbb{R})$ such that $F(e) = V$, that is, $f(T^n(e)) = V_n$ for every $n \in \mathbb{Z}$. □

### 1.5.2 The Hull of a Limit-Periodic Potential

Given a limit-periodic sequence $V \in \ell^\infty(\mathbb{Z})$ with an infinite frequency integer set $S_V \equiv \{n_j\}$, we see that $\text{hull}(V)$ is a Cantor group admitting a minimal translation. Consider the directed set $I = \mathbb{Z}^+$ with the usual order. This gives rise to an inverse system $(\mathbb{Z}_{n_i}, \pi_{ij})_{j\geq i}$, where $\mathbb{Z}_{n_j}$ are $n_j$-cyclic groups with the discrete topology and $\pi_{ij}$ is a homomorphism defined by $\pi_{ij}(k + n_j\mathbb{Z}) = k + n_i\mathbb{Z},\ k \in \mathbb{Z}$.

We endow $\mathbb{Z}_{n_j}$ with a discrete metric defined by $\text{dist}_j(a_1, a_2) = 0$ when $a_1 = a_2$ and $\text{dist}_j(a_1, a_2) = 1$ when $a_1 \neq a_2$. Consider the product group

$$A \equiv \prod_{j=1}^{\infty} \mathbb{Z}_{n_j},$$

of which the topology is the product topology. We endow $A$ with a metric defined by

$$\text{dist}(x, y) = \sum_{j=1}^{\infty} \frac{1}{2^j} \frac{\text{dist}_j(x_j, y_j)}{1 + \text{dist}_j(x_j, y_j)}, \quad x,\ y \in A, \tag{1.6}$$

which is compatible with the product topology.

Let $E \equiv (1, 1, \cdots, 1, \cdots) \in A$, and consider the closed subgroup

$$\bar{B} \equiv \overline{\{nE = (n, n, \cdots, n, \cdots) \in A : n \in \mathbb{Z}\}}.$$

Obviously, $\bar{B}$ is a Cantor group with a minimal translation $T(x) = x + E$, and $\vec{0} \in \bar{B}$ is the identity element. Let

$$\tilde{V}_k \equiv \text{dist}(kE, \vec{0}) = \sum_{j=1}^{\infty} \frac{1}{2^j} \frac{\text{dist}_j(k, 0)}{1 + \text{dist}_j(k, 0)}.$$

By the proof of Lemma 1.33, we know that $\tilde{V} \equiv (\tilde{V}_k)_{k \in \mathbb{Z}}$ is limit-periodic and

$$\text{hull}(\tilde{V}) \cong \bar{B}.$$

Let $P_j(k) \equiv \frac{1}{2^j} \frac{\text{dist}_j(k,0)}{1+\text{dist}_j(k,0)}$. Then we have

$$\tilde{V}(k) = \sum_{j=1}^{\infty} P_j(k),$$

which tells us that one frequency integer set of $\tilde{V}$ is $S_V = \{n_j\}$. By Theorem 1.22, we conclude that

$$\text{hull}(V) \cong \text{hull}(\tilde{V}) \cong \bar{B}.$$

Let $\bar{B}_k = \overline{\{nn_kE : n \in \mathbb{Z}\}} \subset \bar{B}$, and it is easy to see that there exists a decreasing sequence of Cantor subgroups $\bar{B}_k$ with the index $n_k$ and $\bigcap \bar{B}_k = \{\vec{0}\}$.

**Proposition 1.36.** *$b = (k, k, \cdots, k, \cdots)$ is a generator in $\bar{B}$, that is, $\{nb : n \in \mathbb{Z}\}$ is dense in $\bar{B}$, if and only if for any $n_j$, $k$ and $n_j$ have no common divisors.*

*Proof.* If there exists some $n_t$ such that $(k, n_t) = k_t > 1$, then $k$ cannot be a generator for $\mathbb{Z}_{n_t}$, i.e. $nk \neq 1(\text{mod } n_t)$ for any $n \in \mathbb{Z}$, and then

$$\inf_{n\in\mathbb{Z}} \parallel nb - E \parallel \geq \frac{1}{2^t}\frac{\text{dist}_t(nk,1)}{1+\text{dist}_t(nk,1)} = \frac{1}{2^{t+1}},$$

where $E = (1, 1, \cdots, 1, \cdots)$. So $b$ is not a generator.

Conversely, if for any $n_j$, $k$ and $n_j$ have no common divisors, i.e. $(k, n_j) = 1$, we will show that there exists a positive integer sequence $\{q_j\}_{j\in\mathbb{Z}^+}$ such that $\lim_{j\to\infty} q_j b = E$ in the norm sense. Since $(k, n_1) = 1$, there exists some positive integer $q_1$ such that $q_1 k = 1 \pmod{n_1}$. Similarly, since $(k, n_2) = 1$, we have $q_2 k = 1 \pmod{n_2}$. Since $n_1 | n_2$, we still have $q_2 k = 1 \pmod{n_1}$. By induction, we get a sequence $\{q_j\}_{j\in\mathbb{Z}^+}$ such that $q_j k = 1 \pmod{n_i}$ when $i \leq j$. It is easy to see that $\lim_{j\to\infty} q_j b = E$ in the norm sense, and so $b$ is a generator. □

**Proposition 1.37.** *Let $T : \bar{B} \to \bar{B}$, $x \to b + x$. $T$ is minimal if and only if $b$ is a generator.*

*Proof.* If $b$ is not a generator, clearly $T$ is not minimal. On the other hand, if $b$ is a generator, for each $x \in \bar{B}$, there exists a corresponding sequence $\{q_j\}_{j\in\mathbb{Z}^+}$ such that $\lim_{j\to\infty} q_j b = E - x$ ($E$ is the same as before), and so $\lim_{j\to\infty}(q_j b + x) = E$, which implies $\lim_{j\to\infty} T^{q_j}(x) = E$. It follows that $T$ is minimal. □

Furthermore, since the following diagram is commutative where $\pi_i$ is the $i$-th coordinate projection, $(\bar{B}, \pi_i)$ is compatible with this inverse system $(\mathbb{Z}_{n_i}, \pi_{ij})_{j\geq i}$. Proposition 1.27 ensures that $(\bar{B}, \pi_i)$ is also the inverse limit of this system (see the statement (3) of Proposition 1.27). As already introduced, we call $\bar{B}$ a procyclic group, that is, an inverse limit of finite cyclic groups. Equivalently, a procyclic group is a profinite group that can be generated by one element. Thus, we have proved that *a Cantor group with a minimal translation is a procyclic group.* Also, equivalently, we have the following.

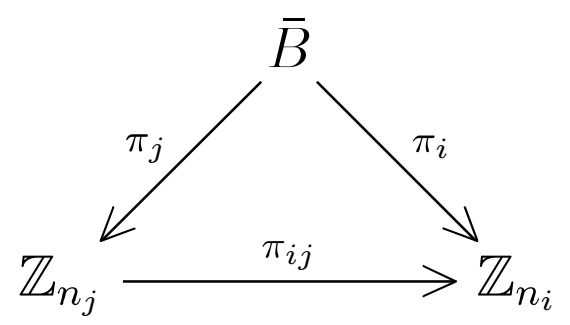

**Theorem 1.38.** *For any limit-periodic potential $V \in \ell^\infty(\mathbb{Z})$,* hull$(V)$ *is a procyclic group.*

By the above discussion and Proposition 1.29, Theorem 1.22 follows directly. There is also a classification theorem about procyclic groups, which is essentially the same as Theorem 1.23. Let $n = \prod_p p^{n(p)}$ be a supernatural number where $p$ goes through all the primes with $0 \leq n(p) \leq \infty$ (one can consider this expression as the generalized prime factorization, and we call it *supernatural* since it is an extension of natural numbers to infinities that differ by different factorizations). We have

**Theorem 1.39.** [38, Theorem 2.7.2] *There exists a unique procyclic group $G$ of order $n$ up to isomorphism.*

**Remark 1.40.** *(1). If $n = p^{n(p)}$, where $p$ is a prime and $n(p) = \infty$, then the associated procyclic group is the group of $p$-adic integers (see p.26 of [38]).*
*(2). Given a limit-periodic potential $V \in \ell^\infty(\mathbb{Z})$,* hull$(V)$ *has order $n$. If $S_V \equiv \{n_j\}$ is a frequency integer set of $V$, then one must have $\lim_{j\to\infty} n_j = n$ (here "=" means that they have the same generalized prime factorization).*

At the end of this section, we discuss quotient groups of a procyclic group. The class of cyclic groups is closed with respect to quotients, that is, a quotient group of a cyclic group is still cyclic. By Proposition 1.31, we know that a quotient group of a procyclic group is still procyclic. We have

**Proposition 1.41.** *Given a procyclic group $G$ with order $n = \prod_{j\in\mathbb{Z}^+} p_j^{r_j}$ where $p_j$ are primes with $0 < r_j \leq \infty$, quotient groups of $G$ are procyclic groups with order $m = \prod_{j\in\mathbb{Z}^+} p_j^{a_j}$ where $0 \leq a_j \leq r_j$, and vice versa.*

*Proof.* Obviously, $G \cong \lim_{\leftarrow k} \mathbb{Z}_{n_k}$, where $n_k | n_{k+1}$ and $\lim_{k\to\infty} n_k = n$. Write $\bar{B} \equiv \overline{\{nE : n \in \mathbb{Z}\}}$, where $E \equiv (1, 1, \cdots, 1, \cdots) \in \prod_k \mathbb{Z}_{n_k}$. Then, $\bar{B} \cong \lim_{\leftarrow k} Z_{n_k}$. It is sufficient to consider $\bar{B}$.

Write $\bar{B}_k \equiv \overline{\{nn_kE : n \in \mathbb{Z}\}}$. Clearly, $\{\bar{B}_k\}$ is a decreasing sequence of open (also closed) subgroups of $\bar{B}$ with index $n_k$ and $\bigcap \bar{B}_k = \{\vec{0}\}$. By Proposition 1.31, for every closed subgroup $N \subset \bar{B}$ (every subgroup is normal in an Abelian group), we have $\bar{B}/N \cong \lim_{\leftarrow k} \bar{B}/\bar{B}_k N$. Since $\bar{B}/\bar{B}_k$ is an $n_k$-cyclic group, $\bar{B}/\bar{B}_k N$ is an $m_k$-cyclic group with $m_k|n_k$. $B/N$ is a procyclic group with order $m = \lim_{k\to\infty} m_k$. Since $m_k|n_k$, it follows that $m = \prod_{j\in\mathbb{Z}^+} p_j^{a_j}$ where $0 \le a_j \le r_j$.

Conversely, if $m = \prod_{j\in\mathbb{Z}^+} p_j^{a_j}$ where $0 \le a_j \le r_j$, obviously there exist $m_k$ such that $m_k|n_k$ and $\lim_{k\to\infty} m_k = m$. Let $\tilde{B} \equiv \lim_{\leftarrow k} \mathbb{Z}_{m_k}$, and define $\phi : \bar{B} \to \tilde{B}$, $\phi(E) = \tilde{E}$, where $\tilde{E} \equiv (1, 1, \cdots, 1, \cdots) \in \prod_k \mathbb{Z}_{m_k}$. Metrics on $\tilde{B}$ and $\bar{B}$ are the metric like (1.6). It is easy to see that $\phi$ is a continuous surjective homomorphism. So $\tilde{B}$ is a quotient group of $\bar{B}$. □

**Remark 1.42.** *By the above proposition, it is easy to see that there exists a universal procyclic group with order $n = \prod_p p^\infty$, where $p$ goes though all the primes, in the sense that any procyclic group is a quotient group of this group.*

# Chapter 2

# Spectral Properties of Limit-Periodic Schrödinger Operators

Examples with purely absolutely continuous spectrum have been known for a long time, dating back to works of Avron and Simon [3], Chulaevsky [6], and Pastur and Tkachenko [33, 34] in the 1980s. There were also works like [29, 30] focusing on pure point spectrum in the class of limit-periodic Schrödinger operators. But, compared with quasi-periodic Schrödinger operators, limit-periodic ones didn't draw that much attention from the research community. In 2008, Avila uploaded a preprint [1] into arXiv (it later appeared in Commun. Math. Phys.), in which he brought the concept of Cantor group into the study of limit-periodic Schrödinger operators. Following Avila's idea, we've exploited the connection between limit-periodic sequences and Cantor groups in the last chapter. In this chapter, we will use the word "procyclic" instead of "Cantor", and utilize such a connection to present properties of limit-periodic Schrödinger operators. Moreover, we will answer questions like, "how often this type of spectrum will occur?" For this goal, the Baire category theorem is helpful; that is, every complete metric space is a Baire space in which for each countable collection of open dense sets $U_n$, their intersection $\bigcap U_n$ is still dense.

Limit-periodic potentials [1] are uniform limits of periodic potentials, so periodic Schrödinger operators play a central role in understanding limit-periodic Schrödinger operators. By Floquet's theorem (also known as Bloch's theorem), the properties of periodic Schrödinger operators are well known. We will first introduce properties of periodic Schrödinger operators, and in the subsequent sections we will show properties about limit-periodic ones, that is the main goal of this thesis.

[1] In this chapter we will use the word "potential" instead of "sequence".

## 2.1 Properties of Periodic Schrödinger Operators

In this section we will focus on the properties of periodic operators that are of interest to us. Please refer to [27, 42, 44, 43] for more details.

Consider a potential $V : \mathbb{Z} \to \mathbb{R}$ that is periodic with period $p$. This gives rise to a periodic Schrödinger operator in $\ell^2(\mathbb{Z})$, given by

$$(H\psi)(n) = \psi(n+1) + \psi(n-1) + V(n)\psi(n). \tag{2.1}$$

We will link the spectral properties of $H$ to properties of the solutions of the associated difference equation

$$u(n+1) + u(n-1) + V(n)u(n) = Eu(n). \tag{2.2}$$

For $k \in \mathbb{R}$ and $l \in \mathbb{Z}$, we define

$$J_l(k) = \begin{pmatrix} V(l) & 1 & & & & e^{-ikp} \\ 1 & V(l+1) & 1 & & & \\ & 1 & V(l+2) & 1 & & \\ & & \ddots & \ddots & \ddots & \\ & & & & & 1 \\ e^{ikp} & & & & 1 & V(l+p-1) \end{pmatrix}.$$

These matrices are the restrictions of $H$ to intervals of length $p$ with suitable self-adjoint boundary conditions. The importance of this choice of boundary condition lies in its connection to the existence of special solutions of (2.2).

The following proposition summarizes several important results concerning the operator, the difference equation, and the matrices. The results listed here are usually referred to as Floquet-Bloch theory.

**Proposition 2.1.** (i). *We have* $E \in \sigma(H)$ *if and only if* (2.2) *has a solution* $\{u(n)\}$ *obeying*

$$u(n+p) = e^{ikp}u(n)$$

*for all* $n$ *and some real number* $k$. *In this case,* $\tilde{u} = \langle u(n) \rangle_{n=l}^{l+p-1}$ *is an eigenvector of the matrix* $J_l(k)$ *corresponding to eigenvalue* $E$.
(ii). *The* $p$ *eigenvalues of* $J_l(k)$ *are independent of* $l$ *and*

$$\sigma(H) = \bigcup_k \sigma(J_l(k)).$$

(iii). *The characteristic polynomial of* $J_l(k)$ *obeys*

$$\det(E - J_l(k)) = \Delta(E) - 2\cos kp,$$

*where* $\Delta(E) = \operatorname{Tr} A_p^{(E,V)}$, *where* $A_p^{(E,V)}$ *is the transfer matrix. We have*

$$\sigma(H) = \{E : |\Delta(E)| \leq 2\}.$$

*The set* $\sigma(H)$ *is made of* $p$ *bands such that on each band,* $\Delta(E)$ *is either strictly increasing or strictly decreasing.*
(iv). *If* $E$ *is in the boundary of some band, we have* $\Delta(E) = \pm 2$. *Moreover, if two different bands intersect, then their common boundary point satisfies* $A_p^{(E,V)} = \pm I$.

*Proof.* For part (i), please see [28] for proof. For parts (ii)—(iv), please see [27, Proof of Proposition 2.1]. □

The function $\Delta$ is called the discriminant associated with the periodic potential $V$. It is a polynomial of degree $p$ with real coefficients.

The other important consequence of periodicity is the existence of a direct integral decomposition. This will be described next. Please refer to [42, Section 5.3] for complete details.

As we have seen above, we can treat $E \in \sigma(H)$ as a function of the variable $k \in [0, \frac{\pi}{p}]$. For each band, the association $k \mapsto E$ is one-to-one and onto. Moreover, if we consider energies in the interior of a band, that is, with $\Delta(E) \in (-2, 2)$ or $k \in (0, \frac{\pi}{p})$, then there are linearly independent solutions $\varphi^{\pm}(E)$ of (2.2) with

$$\varphi^{\pm}_{n+lp}(E) = e^{\pm ilkp}\varphi^{\pm}_n(E).$$

It is easy to see that one can normalize these solutions by requiring

$$\varphi_0^\pm(E) > 0$$

and

$$\sum_{j=0}^{p-1} \left|\varphi_j^\pm(E)\right|^2 = 1. \tag{2.3}$$

With this normalization, we have

$$\varphi^-(E) = \overline{\varphi^+(E)}.$$

Next, we define for $u = \{u_n\}_{n\in\mathbb{Z}}$ of finite support,

$$\hat{u}^\pm(E) = \sum_{n\in\mathbb{Z}} \overline{\varphi_n^\pm(E)} u_n.$$

We also define the measure $d\rho$ on $\sigma(H)$ by

$$d\rho(E) = \frac{1}{\pi}\left|\frac{dk}{dE}(E)\right| \, dE.$$

Then, we have the following result.

**Proposition 2.2.** *The map $u \mapsto \hat{u}$ extends to a unitary map from $\ell^2(\mathbb{Z})$ to $L^2(\sigma(H), d\rho; \mathbb{C}^2)$. Its inverse is given by*

$$\left(\check{f}\right)_n = \frac{1}{2}\int_{\sigma(H)} \left[\varphi_n^+(E) f^+(E) + \varphi_n^-(E) f^-(E)\right] \, d\rho(E).$$

*Moreover, we have that*

$$\widehat{Hu}^\pm(E) = E\hat{u}^\pm(E).$$

*Proof.* This is [42, Theorem 5.3.8.]. □

We use $f^\pm(E)$ for the two components of a $\mathbb{C}^2$-valued function $f \in L^2(\sigma(H), d\rho; \mathbb{C}^2)$.

Proposition 2.2 shows that $H$ has purely absolutely continuous spectrum (of multiplicity two). More precisely, the spectral measure associated with the operator $H$ with periodic potential $V$ and a finitely supported $u \in \ell^2(\mathbb{Z})$ is given by

$$d\mu_{V,u}(E) = g_{V,u}(E)\,dE \tag{2.4}$$

with density

$$g_{V,u}(E) = \frac{1}{2\pi}\left(|\hat{u}^+(E)|^2 + |\hat{u}^-(E)|^2\right)\left|\frac{dk}{dE}(E)\right| \tag{2.5}$$

for $E \in \sigma(H)$ (we set $g_{V,u}(E)$ equal to zero outside of $\sigma(H)$).

We now derive some consequences of the properties of periodic Schrödinger operators which we will use in the next section.

**Lemma 2.3.** *For every $t \in (1,2)$, there exists a constant $D = D(\|V\|_\infty, p, t)$ such that*

$$\int_{\sigma(H)} \left|\frac{dk}{dE}(E)\right|^t \, dE \leq D. \tag{2.6}$$

*Proof.* By Proposition 2.1, we have

$$\left|\frac{dk}{dE}(E)\right| = \left|\frac{\Delta'(E)}{2p\,\sin(kp)}\right|.$$

Since we can bound $|\Delta'(E)|$ by a $(\|V\|_\infty, p)$-dependent constant and $\int_0^\pi (\sin(x))^{1-t}\,dx < \infty$, we have the following estimates,

$$\int_{\sigma(H)} \left|\frac{dk}{dE}(E)\right|^t \, dE \lesssim \int_0^{\frac{\pi}{p}} \left|\frac{1}{2p\,\sin(kp)}\right|^{t-1} dk \lesssim \int_0^{\frac{\pi}{p}} |\sin(kp)|^{1-t}\,dk$$

and the last integral may be bounded by a $t$-dependent constant. ☐

**Lemma 2.4.** *Let $u \in \ell^2(\mathbb{Z})$ have finite support. Then, for every $t \in (1,2)$, there exists a constant $Q = Q(\|V\|_\infty, p, u, t)$ such that*

$$\int_{\sigma(H)} |g_{V,u}(E)|^t \, dE \leq Q. \tag{2.7}$$

*Proof.* Since $u$ has a finite support, we can find a constant $M = M(p, u)$ such that $|\hat{u}^{\pm}(E)|^2 \leq M$. Thus, by (2.5) we have

$$\begin{aligned}\int_{\sigma(H)} |g_{V,u}(E)|^t \, dE &= \int_{\sigma(H)} \left( \frac{1}{2\pi} \left( |\hat{u}^+(E)|^2 + |\hat{u}^-(E)|^2 \right) \left| \frac{dk}{dE}(E) \right| \right)^t dE \\ &\leq \left( \frac{M}{\pi} \right)^t \int_{\sigma(H)} \left| \frac{dk}{dE}(E) \right|^t dE \\ &\leq \left( \frac{M}{\pi} \right)^t D\end{aligned}$$

with the constant $D$ from Lemma 2.3. ☐

**Lemma 2.5.** *Let $(X, d\mu)$ be a finite measure space, let $r > 1$ and let $f_n, f \in L^r$ with $\sup_n \|f_n\|_r < \infty$. Suppose that $f_n(x) \to f(x)$ pointwise almost everywhere. Then, $\|f_n - f\|_p \to 0$ for every $p < r$.*

*Proof.* This is [3, Lemma 2.6]. ☐

**Lemma 2.6.** *Suppose $u \in \ell^2(\mathbb{Z})$ has finite support and $V_n, V : \mathbb{Z} \to \mathbb{R}$ are $p$-periodic and such that $\|V_n - V\|_\infty \to 0$ as $n \to \infty$. Then, for any $t \in (1, 2)$, we have*

$$\int_{\mathbb{R}} |g_{V_n,u}(E) - g_{V,u}(E)|^t \, dE \to 0$$

*as $n \to \infty$.*

*Proof.* By Lemmas 2.4 and 2.5 we only need to prove pointwise convergence. Given the explicit identity (2.5), pointwise convergence follows readily from the following two facts: the discriminant of the approximants converges pointwise to the discriminant of the limit and the matrices $A_l^k$ associated with the approximants converge pointwise to those associated with the limit and therefore so do the associated eigenvectors. ☐

## 2.2 Absolutely Continuous Spectrum

In this section we discuss the absolutely continuous spectrum of limit-periodic Schrödinger operators. We establish the following theorem, which is a kind of discrete version of results established for continuum Schrödinger operators in the 1980's; compare [3, 7, 28].

**Theorem 2.7** (Absolutely Continuous Spectrum)**.** *Suppose $\Omega$ is a procyclic group and $T : \Omega \to \Omega$ a minimal translation. For a dense set of $f \in C(\Omega, \mathbb{R})$ and every $\omega \in \Omega$, the spectrum of $H_\omega$ given by (1.2) is a Cantor set of positive Lebesgue measure and $H_\omega$ has purely absolutely continuous spectrum.*

**Remark 2.8.** *Theorem 2.7 has been published as [12, Theorem 1.1 (a)].*

*A Cantor set is by definition a closed set with empty interior and no isolated points.* For ergodic Schrödinger operators, we've already known that the spectrum is closed and has no isolated points, so the only property that needs to be addressed in the proof of Cantor set is the empty interior.

If instead of $f$ one also considers the one-parameter family $\{\lambda f\}_{\lambda>0}$, then the same conclusion holds for the family uniformly in $\lambda$. This is inferred easily from the proof. In our context we will consider periodic potentials generated from a procyclic group.

**Definition 2.9.** Suppose $\Omega$ is a procyclic group and $T : \Omega \to \Omega$ a minimal translation. We say that a sampling function $f \in C(\Omega, \mathbb{R})$ is $n$-periodic with respect to $T$ if $f(T^n(\omega)) = f(\omega)$ for every $\omega \in \Omega$.

**Proposition 2.10.** *Let $f \in C(\Omega, \mathbb{R})$. If $f(T^{n_0+m}(\omega_0)) = f(T^m(\omega_0))$ for some $\omega_0 \in \Omega$, some minimal translation $T : \Omega \to \Omega$ and every $m \in \mathbb{Z}$, then for every minimal translation $\tilde{T} : \Omega \to \Omega$, $f$ is $n_0$-periodic with respect to $\tilde{T}$.*

*Proof.* Let $\varphi : \Omega \to \ell^\infty(\mathbb{Z})$, $\varphi(\omega) = (f(T^n(\omega)))_{n\in\mathbb{Z}}$. Since $T$ is minimal, the closure of $\{T^n(\omega_0) : n \in \mathbb{Z}\}$ is $\Omega$. By Lemma 1.32 we have $\varphi(\Omega) = \text{hull}(\varphi(\omega_0))$. Since $f(T^{n_0+m}(\omega_0)) = f(T^m(\omega_0))$ for any $m \in \mathbb{Z}$, $\text{hull}(\varphi(\omega_0))$ is a finite set. Then for any $\omega \in \Omega$, $(f(T^n(\omega)))_{n\in\mathbb{Z}}$ is some element in $\text{hull}(\varphi(\omega_0))$. Since every element in $\text{hull}(\varphi(\omega_0))$ is $n_0$-periodic, $(f(T^n(\omega)))_{n\in\mathbb{Z}}$ is $n_0$-periodic. This shows that $f$ is $n_0$-periodic with respect to $T$. That is, we have $f(T^{n_0+m}(\omega)) = f(T^m(\omega))$ for every $\omega \in \Omega$ and $m \in \mathbb{Z}$.

Assume $T$ is the minimal translation by $\omega_1$ and let $\tilde{T}$ be another minimal translation by $\omega_2$. By the previous analysis, we have $f(\omega_1^{n_0+m} \cdot \omega) = f(\omega_1^m \cdot \omega)$ for every $m \in \mathbb{Z}$ and every $\omega \in \Omega$. If $\omega_2$ is equal to $\omega_1^q$ for some integer $q$, obviously we have $f(\tilde{T}^{n_0}(\omega)) = f((\omega_1^q)^{n_0} \cdot \omega) = f(\omega)$ for any $\omega \in \Omega$. If not, since $\{\omega_1^n : n \in \mathbb{Z}\}$ is dense in $\Omega$ (this follows from the minimality of $T$), we have $\lim_{k\to\infty} \omega_1^{n_k} = \omega_2$, and then $f(\omega_2^{n_0} \cdot \omega) = \lim_{k\to\infty} f((\omega_1^{n_k})^{n_0} \cdot \omega) = f(\omega)$. The result follows. $\square$

The above proposition tells us that the periodicity of $f$ is independent of $T$, so we can say $f$ is $n$-periodic without making a minimal translation explicit.

**Proposition 2.11.** *Let $f$ be $p$-periodic. Then, for every $\omega \in \Omega$,*

$$\begin{aligned} L(E,T,f) &= \lim_{m\to\infty} \frac{1}{m} \log \|A_m^{(E,T,f)}(\omega)\| \\ &= \frac{1}{p} \log \rho(A_p^{(E,T,f)}(e)), \end{aligned} \tag{2.8}$$

*where $\rho(A_p^{(E,T,f)}(e))$ is the spectral radius of $A_p^{(E,T,f)}(e)$. In particular, if restricted to periodic sampling functions, the Lyapunov exponent is a continuous function of both the energy $E$ and the sampling function.*

*Proof.* If $f$ is $p$-periodic, it is easy to see that for every $\omega$, $(f(T^n(\omega)))_{n\in\mathbb{Z}}$ is some element of the orbit of $(f(T^n(e)))_{n\in\mathbb{Z}}$, and so its monodromy matrix (i.e., the transfer matrix over one period) is a cyclic permutation of the monodromy matrix associated with $f(T^n(e))$. Thus $\mathrm{Tr} A_p^{(E,T,f)}(\omega)$ is independent of $\omega$, and since $\det A_p^{(E,T,f)}(\omega) = 1$, we can conclude that the eigenvalues of $A_p^{(E,T,f)}(\omega)$ are independent of $\omega$. So the logarithm of the spectral radius of $A_p^{(E,T,f)}(\omega)$ is independent of $\omega$ and (2.8) follows. The continuity statement follows readily. □

We have the following lemma.

**Lemma 2.12.** *Let $f_n \in C(\Omega,\mathbb{R})$ be a sequence of periodic sampling functions converging uniformly to $f_\infty \in C(\Omega,\mathbb{R})$. Assume $\lim_{n\to\infty} L(E,T,f_n)$ exists for every $E$ and the convergence is uniform. Then we have that $L(E,T,f_\infty)$ coincides with $\lim_{n\to\infty} L(E,T,f_n)$ everywhere.*

*Proof.* Since $\lim_{n\to\infty} L(E,T,f_n)$ exists everywhere, from [1, Lemma 2.5], we have $L(E,T,f_n) \to L(E,T,f_\infty)$ in $L^1_{loc}$. So $L(E,T,f_\infty)$ coincides with $\lim_{n\to\infty} L(E,T,f_n)$ almost everywhere. From Proposition 2.11, $L(E,T,f_n)$ is a continuous function, and by uniform convergence, $\lim_{n\to\infty} L(E,T,f_n)$ is also a continuous function. Since $L(E,T,f_\infty)$ is a subharmonic function (cf. [10, Theorem 2.1]), we get that $L(E,T,f_\infty) = \lim_{n\to\infty} L(E,T,f_n)$ for every $E$. The statement follows. □

Next we recall from [1] how periodic sampling functions in $C(\Omega, \mathbb{R})$ can be constructed. Given a procyclic group $\Omega$, a compact subgroup $\Omega_0$ with finite index, and $f \in C(\Omega, \mathbb{R})$, we can define a periodic $f_{\Omega_0} \in C(\Omega, \mathbb{R})$ by

$$f_{\Omega_0}(\omega) = \int_{\Omega_0} f(\omega \cdot \tilde{\omega}) \, d\mu_{\Omega_0}(\tilde{\omega}).$$

Here, $\mu_{\Omega_0}$ denotes Haar measure on $\Omega_0$. This shows that the set of periodic sampling functions is dense in $C(\Omega, \mathbb{R})$. Moreover, as already noted in [1], there exists a decreasing sequence of procyclic subgroups $\Omega_k$ with finite index $n_k$ such that $\bigcap \Omega_k = \{e\}$, where $e$ is the identity element of $\Omega$. (This point has been mentioned in Subsection 1.5.2. That is, for a procyclic group $\bar{B} \equiv \overline{\{nE : n \in \mathbb{Z}\}}$ one can construct subgroups by $\bar{B}_k \equiv \overline{\{nn_kE : n \in \mathbb{Z}\}}$. Clearly, $\{\bar{B}_k\}$ is a decreasing sequence of closed subgroups of $\bar{B}$ with index $n_k$ and $\bigcap \bar{B}_k = \{E\}$). Let $P_k$ be the set of sampling functions defined on $\Omega/\Omega_k$, that is, the elements in $P_k$ are $n_k$-periodic potentials. Denote by $P$ the set of all periodic sampling functions. Then, we have $P_k \subset P_{k+1}$ (which implies $n_k \mid n_{k+1}$) and $P = \bigcup P_k$.

We have the following theorem about the Cantor spectrum, which has been published as [12, Theorem 4.1].

**Theorem 2.13** (Cantor Spectrum)**.** *Suppose $\Omega$ is a procyclic group and $T : \Omega \to \Omega$ a minimal translation. There exists a dense $G_\delta$ set $\mathcal{C} \subseteq C(\Omega, \mathbb{R})$ such that for every $f \in \mathcal{C}$ and $\omega \in \Omega$, the spectrum of the operator $H_\omega$ given by (1.2) is a Cantor set.*

*A $G_\delta$ set is a countable intersection of open subsets.* Fix $\Omega$ and $T$ as in the theorem throughout this section. By minimality, for given $f \in C(\Omega, \mathbb{R})$, the spectrum of $H_\omega$ is independent of $\omega$. For notational convenience, we will denote the spectrum by $\Sigma(f)$.

**Lemma 2.14.** $\mathfrak{N} := \{f \in C(\Omega, \mathbb{R}) : \Sigma(f)$ *has empty interior*$\}$ *is a $G_\delta$ set.*

*Proof.* This is essentially [3, Lemma 1.1]. □

**Lemma 2.15.** *For every $f \in P_k$, $k \in \mathbb{N}$ and every $\varepsilon > 0$, there exists $\tilde{f}$ in $P_k$ satisfying $\|f - \tilde{f}\| < \varepsilon$ such that $\Sigma(\tilde{f})$ has exactly $p_k$ components, that is, its $p_k - 1$ gaps are all open.*

*Proof.* This follows from the proof of [1, Claim 3.4]. For the reader's convenience, we provide a proof. Let $f \in P_k$ be given. By $\omega$-independence of the spectrum, we may choose and fix an arbitrary $\omega \in \Omega$ for the purpose of this proof. Next, given $\varepsilon > 0$, let $M$ be large enough so that $\frac{2p_k+1}{M} < \varepsilon$. Then, for $1 \leq t \leq 2p_k + 1$, there is $\tilde{f}_t \in P_k$ with

$$\tilde{f}_t(T^i\omega) = f(T^i\omega),\ 0 \leq\ i\ \leq p_k - 2 \ \text{ and }\ \tilde{f}_t(T^{p_k-1}\omega) = f(T^{p_k-1}\omega) + \frac{t}{M}.$$

Obviously, $\|\tilde{f}_t - f\| < \varepsilon$ and we claim that there exists some $t$ in this range such that the spectrum associated with $\tilde{f}_t$ has exactly $p_k$ components.

Suppose this claim fails. Then, for every $t$ in this range, there exists by Proposition 2.1.(d) an energy $E_t \in \Sigma(\tilde{f}_t)$ with $T_{p_k}^{(E_t,\tilde{f}_t)} = \pm id$. That is,

$$\begin{pmatrix} E_t - f(T^{p_k-1}\omega) - \frac{t}{M} & -1 \\ 1 & 0 \end{pmatrix} \begin{pmatrix} E_t - f(T^{p_k-2}\omega) & -1 \\ 1 & 0 \end{pmatrix} \cdots \begin{pmatrix} E_t - f(\omega) & -1 \\ 1 & 0 \end{pmatrix} = \pm id.$$

Since $f$ is $p_k$-periodic, we get from this

$$T_{p_k}^{(E_t,f)} = \pm \begin{pmatrix} 1 & \frac{t}{M} \\ 0 & 1 \end{pmatrix}. \tag{2.9}$$

Indeed, rewriting the identity above, we find that

$$\begin{aligned} T_{p_k}^{(E_t,f)} &= \pm id + \begin{pmatrix} \frac{t}{M} & 0 \\ 0 & 0 \end{pmatrix} \begin{pmatrix} E_t - f(T^{p_k-2}\omega) & -1 \\ 1 & 0 \end{pmatrix} \cdots \begin{pmatrix} E_t - f(\omega) & -1 \\ 1 & 0 \end{pmatrix} \\ &= \pm id + \begin{pmatrix} \frac{t}{M} & 0 \\ 0 & 0 \end{pmatrix} \begin{pmatrix} E_t - f(T^{p_k-1}\omega) & -1 \\ 1 & 0 \end{pmatrix}^{-1} T_{p_k}^{(E_t,f)} \end{aligned}$$

so that

$$\begin{aligned} T_{p_k}^{(E_t,f)} &= \pm \left( id - \begin{pmatrix} \frac{t}{M} & 0 \\ 0 & 0 \end{pmatrix} \begin{pmatrix} E_t - f(T^{p_k-1}\omega) & -1 \\ 1 & 0 \end{pmatrix}^{-1} \right)^{-1} \\ &= \pm \left( id - \begin{pmatrix} \frac{t}{M} & 0 \\ 0 & 0 \end{pmatrix} \begin{pmatrix} 0 & 1 \\ -1 & E_t - f(T^{p_k-1}\omega) \end{pmatrix} \right)^{-1} \\ &= \pm \left( id - \begin{pmatrix} 0 & \frac{t}{M} \\ 0 & 0 \end{pmatrix} \right)^{-1} \\ &= \pm \begin{pmatrix} 1 & -\frac{t}{M} \\ 0 & 1 \end{pmatrix}^{-1} \\ &= \pm \begin{pmatrix} 1 & \frac{t}{M} \\ 0 & 1 \end{pmatrix}. \end{aligned}$$

The relation (2.9) implies that if $t \neq t'$, we have $T_{p_k}^{(E_t,f)} \neq T_{p_k}^{(E_{t'},f)}$, and therefore $E_t \neq E_{t'}$. But there are at most $2p_k$ values of $E$ for which $tr\ T_{p_k}^{(E,f)} = \pm 2$; contradiction. □

**Lemma 2.16.** *Let $f \in C(\Omega, \mathbb{R})$ be $n$-periodic. The Lebesgue measure of each band of $\Sigma(f)$ is at most $\frac{2\pi}{n}$.*

*Proof.* This result is well known. Please refer to [1, Lemma 2.4 (i)] for the proof. □

*Proof of Theorem 2.13 (Cantor Spectrum).* The proof is close in spirit to the proof of [3, Theorem 1]. By Lemma 2.14, we only need to prove that $\mathfrak{N}$ is dense. Since $P = \bigcup P_k$ is dense in $C(\Omega, \mathbb{R})$, it suffices to show that, given $f \in P$ and $\varepsilon > 0$, there is a sampling $\tilde{f}$ such that $\|\tilde{f} - f\| < 2\varepsilon$ and $\Sigma(\tilde{f})$ is nowhere dense.

So let $f \in P$ and $\varepsilon > 0$ be given. Write $f$ as $f = \sum_{j=0}^{N} a_j W_j$ with $W_j \in P_j$. We construct $s_0 = \sum_{i=0}^{N} a_n^{(0)} W_i$ so that $\|s_0\| < \varepsilon$ and $f_0 = f + s_0$ has all $n_N - 1$ gaps open. (This is possible due to Lemma 2.15.)

Suppose $s_0, s_1, \dots, s_{k-1}$ are picked. Let $\alpha_{k-1}$ be the minimal gap size of $f_{k-1}$ and define $\beta_k = \min\{\alpha_0, \alpha_1, \dots, \alpha_{k-1}\}$. Applying Lemma 2.15, we pick $s_k = \sum_{i=0}^{N+k} a_i^{(k)} W_i$ so that

$$\|s_k\| < \frac{\varepsilon}{2^k}, \tag{2.10}$$

$$\|s_k\| < \frac{1}{3}\frac{\beta_k}{2^k}, \tag{2.11}$$

$$f_k = f + \sum_{j=0}^{k} s_j \text{ has all gaps open.} \tag{2.12}$$

The limit of $f_k$ exists by (2.10), let $\tilde{f} = \lim_{k\to\infty} f_k$. By construction, we have $\|\tilde{f} - f\| < \varepsilon$. We claim that $\Sigma(\tilde{f})$ is nowhere dense; equivalently, its complement is dense.

Given $E \in \Sigma(\tilde{f})$ and $\tilde{\varepsilon} > 0$, we can pick $k$ large enough so that

$$\|\tilde{f} - f_k\| < \frac{\tilde{\varepsilon}}{3}, \tag{2.13}$$

$$\frac{2\pi}{n_{N+k}} < \frac{\tilde{\varepsilon}}{3}, \tag{2.14}$$

$$\frac{\varepsilon}{2^{k-1}} < \frac{\tilde{\varepsilon}}{3}. \tag{2.15}$$

By (2.13), there exists $E' \in \Sigma(f_k)$ such that $|E' - E| < \frac{\tilde{\varepsilon}}{3}$. Moreover, by Lemma 2.16 and (2.14), we can find $\tilde{E}$ in a gap of $\Sigma(f_k)$ such that $|E' - \tilde{E}| < \frac{\tilde{\varepsilon}}{3}$. Write this gap of $\Sigma(f_k)$ that contains $\tilde{E}$ as $(a - \delta, a + \delta)$. By definition, we have $2\delta \geq \beta_{k+1}$. By (2.11),

$$\|\tilde{f} - f_k\| = \left\| \sum_{j=k+1}^{\infty} s_j \right\| < \frac{\beta_{k+1}}{3} \left( \frac{1}{2^{k+1}} + \frac{1}{2^{k+2}} + \cdots \right) \leq \frac{\delta}{3},$$

so we have $(a - \frac{\delta}{3}, a + \frac{\delta}{3}) \cap \Sigma(\tilde{f}) = \emptyset$.

We claim that there exists $\delta' \in [\frac{\delta}{3}, \delta)$ such that $(a - \delta', a + \delta') \cap \Sigma(\tilde{f}) = \emptyset$ and $|\delta' - \delta| < \frac{\varepsilon}{2^k}$. As we saw above, we may arrange for $\delta' \geq \frac{\delta}{3}$. Suppose that it is impossible to find such a $\delta'$ with $|\delta' - \delta| < \frac{\varepsilon}{2^k}$. Then, there will be a point $x \in \Sigma(\tilde{f})$ such that $[x - \frac{\varepsilon}{2^k}, x + \frac{\varepsilon}{2^k}] \subseteq (a - \delta, a + \delta)$. Then we get a contradiction to the already established fact $[x - \frac{\varepsilon}{2^k}, x + \frac{\varepsilon}{2^k}] \cap \Sigma(f_k) = \emptyset$ since (2.10) implies

$$\|\tilde{f} - f_k\| = \left\| \sum_{j=k+1}^{\infty} s_j \right\| < \frac{\varepsilon}{2^k}.$$

We can choose an energy $\hat{E}$ in the gap of $\Sigma(\tilde{f})$ that contains $(a - \delta', a + \delta')$ such that $|\hat{E} - \tilde{E}| \leq \frac{\varepsilon}{2^k} < \frac{\tilde{\varepsilon}}{3}$, where the second inequality follows from (2.15). Moreover, since we also have $|\tilde{E} - E'| < \frac{\tilde{\varepsilon}}{3}$ and $|E' - E| < \frac{\tilde{\varepsilon}}{3}$, it follows that $|\hat{E} - E| < \tilde{\varepsilon}$. This shows that $\mathbb{R} \setminus \Sigma(\tilde{f})$ is dense and completes the proof. □

**Remark 2.17.** *The reader may notice that the statement of Theorem 2.13 is a consequence of [1, Corollary 1.2] (any zero measure set has no interior points). However, the main purpose here is the method of proof presented here, which is direct and flexible enough so that it can be used to also ensure absolutely continuous spectrum. In particular, the Cantor spectra constructed here may have positive Lebesgue measure, whereas the Cantor spectra generated in the proof of [1, Corollary 1.2] always have zero Lebesgue measure.*

It is now ready to prove our main theorem.

*Proof of Theorem 2.7 (Absolutely Continuous Spectrum).* The idea is to modify the construction from the proof of Theorem 2.13. Thus, we will again start with an arbitrarily small ball in $C(\Omega, \mathbb{R})$ and construct a point in this ball for which the associated Schrödinger operator has both Cantor spectrum and purely absolutely continuous spectrum. The presence of absolutely continuous spectrum then also implies that the Lebesgue measure of the spectrum is positive (for an absolute continuous spectrum cannot be supported on any zero measure set).

Fix $t \in (1,2)$ and let $u \in \ell^2(\mathbb{Z})$ have finite support. In going through the construction in the proof of Theorem 2.13, pick $s_k$ so that in addition to the conditions above, we have

$$\left( \int_{-\infty}^{\infty} \left| g_u^{k-1}(E) - g_u^k(E) \right|^t \right)^{\frac{1}{t}} \leq \frac{1}{2^k}, \tag{2.16}$$

where $g_u^k$ is the density of the spectral measure associated with $u$ and the periodic potential $n \mapsto f_k(T^n\omega)$, with the estimate above being uniform in $\omega \in \Omega$. This is possible due to Lemma 2.6.

By Lemma 2.4, there exists $Q(u,t) < \infty$ such that $\int_{\mathbb{R}} \left| g_u^k(E) \right|^t dE \leq Q(u,t)$.

Fix any $\omega \in \Omega$. Let $A$ be a finite union of open sets. If $P_A^k$ is the spectral projection for the potential $n \mapsto f_k(T^n\omega)$ and $P_A$ is the spectral projection for the potential $n \mapsto \tilde{f}(T^n\omega)$, it follows that $\langle u, P_A u\rangle \leq \limsup_{k\to\infty} \langle u, P_A^k u\rangle$ since $\|f_k - \tilde{f}\|_\infty \to 0$ and hence the associated Schrödinger operators converge in norm.

Applying Hölder's inequality, we find

$$\langle u, P_A u\rangle \leq \limsup_{k\to\infty} \int_A g_u^k(E)\, dE \leq Q(u,t)|A|^{\frac{1}{q}},$$

where $\frac{1}{q} + \frac{1}{t} = 1$ and $|\cdot|$ denotes Lebesgue measure. This shows that the spectral measure associated with $u$ and the Schrödinger operator with potential $n \mapsto \tilde{f}(T^n\omega)$ is absolutely continuous with respect to Lebesgue measure. Since this holds for every finitely supported $u$, it follows that this operator has purely absolutely continuous spectrum. □

Note that the denseness cannot be extended to be generic ($G_\delta$-set), for intersection of generic sets is still generic (by the Baire category theorem) and we will see in the following section that there exists a generic set of $f \in C(\Omega, \mathbb{R})$ presenting purely singular continuous spectrum.

## 2.3 Singular Continuous Spectrum

We will here discuss the singular continuous spectrum. The Lyapunov exponents are always zero for absolutely continuous spectrum (see Theorem 1.7). For singular continuous spectrum, the corresponding Lyapunov exponent may be zero or positive. We will discuss the singular continuous spectrum in the regimes of zero and positive Lyapunov exponents respectively.

### 2.3.1 Zero Lyapunov Exponents

**Theorem 2.18** (Singular Continuous Spectrum I). *Suppose $\Omega$ is a procyclic group and $T : \Omega \to \Omega$ a minimal translation. For a dense $G_\delta$-set of $f \in C(\Omega, \mathbb{R})$, every $\omega \in \Omega$ and every $\lambda \neq 0$, $\Sigma(\lambda f)$ has zero Lebesgue measure, $H_\omega$ has purely singular continuous spectrum for every $\omega \in \Omega$, and $E \mapsto L(E, T, \lambda f)$ is zero.*

The above theorem has been published as [12, Theorem 1.2.(b)]. Zero Lebesgue measure spectrum will automatically preclude absolutely continuous spectrum. We still need to find a way to preclude point spectrum. The Gordon lemma is a powerful weapon for us.

**Definition 2.19.** A bounded function $V : \mathbb{Z} \to \mathbb{R}$ is called a Gordon potential if there are positive integers $q_k \to \infty$ such that

$$\max_{1 \leq n \leq q_k} |V(n) - V(n \pm q_k)| \leq k^{-q_k}$$

for every $k \geq 1$.

Clearly, if $V$ is a Gordon potential, then so is $\lambda V$ for every $\lambda \in \mathbb{R}$.

**Lemma 2.20** (Gordon Lemma). *Suppose $V$ is a Gordon potential. The Schrödinger operator $H$ given by* (2.1) *has empty point spectrum.*

This is essentially due to Gordon [21]; see [17] for the modification of the argument necessary to prove the result as stated.

**Lemma 2.21.** *Suppose $\Omega$ is a procyclic group and $T : \Omega \to \Omega$ a minimal translation. Then there exists a dense $G_\delta$-set $\mathcal{G} \subseteq C(\Omega, \mathbb{R})$ such that for every $f \in \mathcal{G}$ and $\omega \in \Omega$, the potential $V_\omega$ is a Gordon potential.*

*Proof.* We know the set of periodic potentials is dense in $C(\Omega, \mathbb{R})$. For $j, k \in \mathbb{N}$, let

$$\mathcal{G}_{j,k} = \left\{ f \in C(\Omega, \mathbb{R}) : \text{ there is a } j\text{-periodic } f_j \text{ such that } \|f - f_j\| < \frac{1}{2}(jk)^{(-jk)} \right\}.$$

Clearly, $\mathcal{G}_{j,k}$ is open. For $k \in \mathbb{N}$, let

$$\mathcal{G}_k = \bigcup_{j=1}^{\infty} \mathcal{G}_{j,k}.$$

The set $\mathcal{G}_k$ is open by construction and dense since it contains all periodic sampling functions. Thus,

$$\mathcal{G} = \bigcap_{k=1}^{\infty} \mathcal{G}_k$$

is a dense $G_\delta$ subset of $C(\Omega, \mathbb{R})$. We claim that for every $f \in \mathcal{G}$ and every $\omega \in \Omega$, the potential $V_\omega$ is a Gordon potential.

Let $f \in \mathcal{G}$ and $\omega \in \Omega$ be given. Since $f \in A_k$ for every $k \in \mathbb{N}$, we can find $j_k$-periodic $f_{j_k}$ satisfying

$$\|f - f_{j_k}\| < \frac{1}{2}(j_k k)^{-j_k k}.$$

Let $q_k = j_k k$, so that $q_k \to \infty$ as $k \to \infty$. Then, we have

$$\begin{aligned}
\max_{1 \le n \le q_k} \|V_\omega(n) - V_\omega(n \pm q_k)\| &= \max_{1 \le n \le q_k} \|f(T^n\omega) - f(T^{n \pm q_k}\omega)\| \\
&= \max_{1 \le n \le j_k k} \|f(T^n\omega) - f_{j_k}(T^n\omega) + f_{j_k}(T^{n \pm j_k k}\omega) - f(T^{n \pm j_k k}\omega)\| \\
&\le \max_{1 \le n \le j_k k} \|f(T^n\omega) - f_{j_k}(T^n\omega)\| + \max_{1 \le n \le j_k k} \|f(T^{n \pm j_k k}\omega) - f_{j_k}(T^{n \pm j_k k}\omega)\| \\
&< \frac{1}{2}(j_k k)^{-j_k k} + \frac{1}{2}(j_k k)^{-j_k k} \\
&\le k^{-j_k k} \\
&= k^{-q_k}.
\end{aligned}$$

It follows that $V_\omega$ is a Gordon potential. □

**Lemma 2.22.** *Suppose $T : \Omega \to \Omega$ is a minimal translation of a procyclic group. For a dense $G_\delta$ set of $f \in C(\Omega, \mathbb{R})$, and for every $\lambda \neq 0$, the Schrödinger operator with potential $\lambda f(T^n\omega)$ has a spectrum of zero Lebesgue measure for every $\omega \in \Omega$ and meanwhile zero Lyapunov exponent.*

*Proof.* This is [1, Corollary 1.2.]. ☐

We can now give the proof of Theorem 2.18.

*Proof of Theorem 2.18.* Since the intersection of two dense $G_\delta$ sets is again a dense $G_\delta$ set and zero-measure spectrum precludes absolutely continuous spectrum, the result follows directly from Theorem 2.21 and Lemma 2.22. ☐

### 2.3.2 Positive Lyapunov Exponents

**Theorem 2.23** (Singular Continuous Spectrum II)**.** *Suppose $\Omega$ is a procyclic group and $T : \Omega \to \Omega$ a minimal translation. Then for a dense set of $f \in C(\Omega, \mathbb{R})$, every $\omega \in \Omega$ and every $\lambda \neq 0$, $\Sigma(\lambda f)$ has zero Hausdorff dimension, $H_\omega$ has purely singular continuous spectrum for every $\omega \in \Omega$, and $E \mapsto L(E, T, \lambda f)$ is positive.*

**Remark 2.24.** *The above theorem has been published as [13, Theorem 1.3.]*

A key part in proving Theorem 2.23 is to establish the following result:

**Theorem 2.25.** *Suppose $\Omega$ is a Cantor group. Then these exists a dense set $\mathcal{F} \subset C(\Omega, \mathbb{R})$ such that for every $f \in \mathcal{F}$, every minimal translation $T : \Omega \to \Omega$, every $\omega \in \Omega$, and every $\lambda \neq 0$, $\lambda f(T^n(\omega))$ is a Gordon potential.*

First, for our relatively restricted purposes, we will simply recall the definition of Hausdorff measures and Hausdorff dimension. We refer the reader to [36] for more information.

**Definition 2.26.** Let $A \subseteq \mathbb{R}$ be a subset. A countable collection of intervals $\{b_n\}_{n=1}^{\infty}$ is called a $\delta$-cover of $A$ if $A \subset \bigcup_{n=1}^{\infty} b_n$ with $|b_n| < \delta$ for all $n$'s. (Here, $|\cdot|$ denotes Lebesgue measure, and we will adopt this notation throughout the section.)

**Definition 2.27.** Let $\alpha \in \mathbb{R}$. For any subset $A \subseteq \mathbb{R}$, the $\alpha$-dimensional Hausdorff measure of $A$ is defined as

$$h^\alpha(A) = \lim_{\delta \to 0} \inf_{\delta\text{-covers}} \sum_{n=1}^{\infty} |b_n|^\alpha. \tag{2.17}$$

The $h^\alpha(A)$ is well defined as an element of $[0,\infty]$ since $\inf_{\delta\text{-covers}} \sum_{n=1}^{\infty} |b_n|^\alpha$ is monotonically increasing as $\delta$ decreases to zero and therefore the limit in (2.17) exists. Restricted to the Borel sets, $h^1$ coincides with Lebesgue measure and $h^0$ is the counting measure. If $\alpha < 0$, we always have $h^\alpha(A) = \infty$ for any $A \neq \emptyset$, while if $\alpha > 1$, $h^\alpha(\mathbb{R}) = 0$.

It is not hard to see that for every $A \subseteq \mathbb{R}$, there is a unique $\alpha \in [0,1]$, called the Hausdorff dimension $\dim_H(A)$ of $A$, such that $h^\beta(A) = \infty$ for every $\beta < \alpha$ and $h^\beta(A) = 0$ for every $\beta > \alpha$. In particular, every $A \subseteq \mathbb{R}$ with $|A| > 0$ must have $\dim_H(A) = 1$.

Let's now move to the proof of Theorem 2.23. For convenience, we write $A_n^{(E,f)}(\omega) = A_n^{(E,f,T)}(\omega)$, $A_n^{(E,f)} = A_n^{(E,f,T)}(e)$, and $L(E,f) = L(E,T,f)$. Since $T : \Omega \to \Omega$ is a minimal translation, the homomorphism $\mathbb{Z} \to \Omega$, $n \to T^n e$ is injective with dense image in $\Omega$, and we can write $f(n) = f(T^n(e))$ without any conflicts.

**Lemma 2.28.** *Let $f \in C(\Omega,\mathbb{R})$ be $n$-periodic. Let $C \geq 1$ be such that for every $E \in \Sigma(f)$, there exist $\omega \in \Omega$ and $k \geq 1$ such that $\|A_k^{(E,T,f)}(\omega)\| \geq C$. Then, $|\Sigma(f)| \leq \frac{4\pi n}{C}$.*

*Proof.* This is just [1, Lemma 2.4 (ii)]. □

The above lemma use the norm of the transfer matrix to estimate the measure of the spectrum, which is helpful in our proof of hausdorff dimension zero. We also need two more lemmas. More precisely, we will make further use of the constructions which play central roles in the proof of these two lemmas.

**Lemma 2.29.** *Let $B$ be an open ball in $C(\Omega,\mathbb{R})$, let $F \subset P \cap B$ be a finite family of sampling functions, and let $0 < \varepsilon < 1$. Then there exists a sequence $F_K \subset P \cap B$ such that*
*(i). $L(E,\lambda F_K) > 0$ whenever $\varepsilon \leq |\lambda| \leq \varepsilon^{-1}$, $E \in \mathbb{R}$,*
*(ii). $L(E,\lambda F_K) \to L(E,\lambda F)$ uniformly on compacts (as functions of $(E,\lambda) \in \mathbb{R}^2$).*

This is [1, Lemma 3.1]. As in [1], we use the notation

$$L(E,\lambda F) = \frac{1}{\#F} \sum_{f \in F} L(E,T,\lambda f),$$

where $F$ is a finite family of sampling functions (with multiplicities!) and $\lambda \in \mathbb{R}$. The proof of this lemma is constructive. We will describe this construction explicitly in the proof of Theorem 2.7.

**Lemma 2.30.** *Let $B$ be an open ball in $C(\Omega, \mathbb{R})$, and let $F \subset P \cap B$ be a finite family of sampling functions. Then for every $N \geq 2$ and $K$ sufficiently large, there exists $F_K \subset P_K \cap B$ such that*
*(i). $L(E, \lambda F_K) \to L(E, \lambda F)$ uniformly on compacts (as functions of $(E, \lambda) \in \mathbb{R}^2$).*
*(ii). The diameter of $F_K$ is at most $n_K^{-N/2}$.*

This lemma is a variation of [1, Lemma 3.2 ]. We will prove this lemma using suitable modifications of Avila's arguments. Some of these modifications, which will later enable us to prove the Gordon property, are not apparent from the statement of the lemma. We will give detailed arguments for the modified parts of the proof and refer the reader to [1] for the parts that are analogous.

*Proof of Lemma 2.30.* Assume that $F = \{f_1, f_2, \ldots, f_m\} \subset C(\Omega, \mathbb{R})$ is a finite family of $n_k$-periodic sampling functions with $n_k \geq 2$, and let $K > k$ be large enough. We construct $F_K^{\vec{t}}$ as follows. Let $n_K = mn_k r + d$, $0 \leq d \leq mn_k - 1$ and $n_k | d$. Let $I_j = [jn_k, (j+1)n_k - 1] \subset \mathbb{Z}$ and let $0 = j_0 < j_1 < \cdots < j_{m-1} < j_m = n_K/n_k$ be a sequence such that $j_{i+1} - j_i = r + 1$ when $0 \leq i < d/n_k$ and $j_{i+1} - j_i = r$ when $d/n_k \leq i \leq m-1$. Define an $n_K$-periodic $f$ as follows. For $0 \leq l \leq n_K - 1$, let $j$ be such that $l \in I_j$ and let $i$ be such that $j_{i-1} \leq j < j_i$ and then let $f(l) = f_i(l)$. Next, for any sequence $\vec{t} = (t_1, t_2, \ldots, t_m)$ with $t_i \in \{0, 1, \ldots, r-1\}$, we define an $n_K$-periodic $f_K^{\vec{t}}$ as follows. If $j = j_i - 1$ for some $1 \leq i < m$, we let $f_K^{\vec{t}}(l) = f(l) + r^{-N} t_i$, and if $j = j_m - 2$, we let $f_K^{\vec{t}}(l) = f(l) + r^{-N} t_m$. Otherwise we let $f_K^{\vec{t}}(l) = f(l)$. Let $F_K^{\vec{t}}$ be the family consisting of all $f_K^{\vec{t}}$'s. The statement (ii) is clear for large $K$. (Note: in [1], Avila's construction is such that if $j = j_i - 1$ for some $1 \leq i \leq m$, then $f_K^{\vec{t}}(l) = f(l) + r^{-20} t_i$; otherwise, $f_K^{\vec{t}}(l) = f(l)$.)

For fixed $E$ and $\lambda$, we let $A_{n_K}^{(E, \lambda f_K^{\vec{t}})} = C^{(t_m, m)} B^{(m)} \ldots C^{(t_1, 1)} B^{(1)}$, where $B^{(i)} = (A_{n_k}^{(E, \lambda f_i)})^{j_i - j_{i-1} - 1}, 1 \leq i \leq m-1$ and $B^{(m)} = (A_{n_k}^{(E, \lambda f_m)})^{j_m - j_{m-1} - 2}$, and $C^{(t_i, i)} = A_{n_k}^{(E - \lambda r^{-N} t_i, \lambda f_i)}, 1 \leq i \leq m-1$ and $C^{(t_m, m)} = A_{n_k}^{(E, \lambda f_m)} A_{n_k}^{(E - \lambda r^{-N} t_m, \lambda f_m)}$. When $E$ and $\lambda$ are in a compact set, the norm of the $C^{(t_i, i)}$-type matrices is bounded as $r$ grows, while the norm of the $B^{(i)}$-type matrices may get large.

Notice that our perturbation here is $r^{-N}t$ (as opposed to Avila's $r^{-20}t$ perturbation in [1, Lemma 3.2 ]), so [1, Claim 3.7] should be replaced by the following version:

"Let $s_j$ be the most contracted direction of $\hat{B}^{(j)}$ and let $u_j$ be the image under $\hat{B}^{(j)}$ of the most expanded direction. Call $\vec{t}$ $j$-nice, $1 \leq j \leq d$, if the angle between $\hat{C}^{(j)}u_j$ and $S_{j+1}$ (less than $\pi$) is at least $r^{-3N}$ with the convention that $j+1=1$ for $j=d$. Let $r$ be sufficiently large, and let $\vec{t}$ be $j$-nice. If $z$ is a non-zero vector making an angle at least $r^{-4N}$ with $s_j$, then $z' = \hat{C}^{(j)}\hat{B}^{(j)}z$ makes an angle at least $r^{-4N}$ with $S_{j+1}$ and $\|z'\| \geq \|\hat{B}^{(j)}\|r^{-5N}\|z\|$."
The proof of [1, Claim 3.7] can be applied to get the above version of the claim with the corresponding quantitative modification. Moreover, we have also made a little shift in the perturbation, so $C^{(t_m,m)} = A_{n_k}^{(E,\lambda f_m)} A_{n_k}^{(E-\lambda r^{-N}t_m,\lambda f_m)}$, while Avila's $C^{(t_m,m)} = A_{n_k}^{(E-\lambda r^{-20}t_m,\lambda f_m)}$. [1, Claim 3.8] still holds, but Avila's proof of [1, Claim 3.8] cannot be applied directly. To this end we prove the following claim:

**Claim 2.31.** *For every $M \in \mathrm{SL}(2,R)$, there are $m_1, m_2 \in (0,\infty)$ with the following property. Suppose $A$ and $B$ are two vectors in $\mathbb{R}^2$, and $\Delta\theta$ is the angle between $A$ and $B$ with $0 < \Delta\theta \leq \pi$. Let $\Delta\tilde{\theta}$ be the angle between $MA$ and $MB$ (again so that $0 < \Delta\tilde{\theta} \leq \pi$). Then, $m_1\Delta\theta \leq \Delta\tilde{\theta} \leq m_2\Delta\theta$.*

*Proof.* By the singular value decomposition (see [41, Theorem 2.5.1]), there exist $O_1$ and $O_2$ in $\mathrm{SO}(2,R)$ such that $M = O_1SO_2$, where $S$ is a diagonal matrix. Since $O_1$ and $O_2$ are rotations on $\mathbb{R}^2$, it is sufficient to consider

$$S = \begin{pmatrix} \mu_1 & 0 \\ 0 & \mu_1^{-1} \end{pmatrix}.$$

Without loss of generality, assume $\mu_1 \geq 1$. Let $A = (a,b)^t$ ($t$ denotes the transpose of vectors) and $B = (c,d)^t$ be two normalized vectors, and let $\theta_A$ and $\theta_B$ be the argument of $A$ and the argument of $B$ respectively. Let $\tilde{A} = SA = (a\mu_1, b/\mu_1)^t$ with the argument $\theta_{\tilde{A}}$ and $\tilde{B} = SB = (c\mu_1, d/\mu_1)^t$ with the argument $\theta_{\tilde{B}}$.

We adopt the following notation for convenience. Let $I$, $II$,$III$,$IV$ denote one of two vectors in the first quadrant (including $\{(x,0) : x \geq 0\}$), the second quadrant (including $\{(0,y) : y > 0\}$), the third quadrant (including $\{(x,0) : x < 0\}$) and the fourth quadrant (including $\{(0,y) : y < 0\}$), respectively. Then $(I,I)$ denotes that both two vectors are in the first quadrant, $(I, II)$ denotes that one vector is in the first quadrant while the other is in the second quadrant, and so on.

We will need the following observation:

$$0 < \theta_1, \theta_2 < \pi/2 \text{ and } \tan\theta_1 \geq \frac{\tan\theta_2}{\mu_1^2} \quad \Rightarrow \quad \theta_1 \geq \frac{\theta_2}{4\mu_1^2}. \tag{2.18}$$

Indeed, since $\tan\theta_1 \geq \frac{1}{\mu_1^2}\tan\theta_2 \geq \frac{1}{2\mu_1^2}\theta_2$ and $0 < \frac{\theta_2}{2\mu_1^2} < 1$, we have

$$\theta_1 \geq \arctan\frac{\theta_2}{2\mu_1^2} = \frac{\theta_2}{2\mu_1^2} - (\frac{\theta_2}{2\mu_1^2})^3/3 + O((\frac{\theta_2}{2\mu_1^2})^5) \geq \frac{\theta_2}{4\mu_1^2}.$$

For the proof of Claim 2.31, we consider two cases.

*Case 1.* $\pi/2 \leq \Delta\theta \leq \pi$. Here $A$ and $B$ cannot be in the same quadrant. Notice that the impact of $S$ on vectors is to move them closer to the $x$-axis and keep them in the same quadrant. Thus, for the subcases $(I, II)$, $(I, III)$, $(II, IV)$ and $(III, IV)$, we can easily conclude that $\Delta\theta/2 \leq \Delta\tilde{\theta} \leq 2\Delta\theta$. There are two subcases left, $(I, IV)$ and $(II, III)$. We will discuss $(I, IV)$; the method can be readily adapted to $(II, III)$. For $(I, IV)$, if $\theta_A = 0$ and $\theta_B = 3\pi/2$, then $\theta_{\tilde{A}}$ and $\theta_{\tilde{B}}$ are also 0 and $3\pi/2$ respectively, and so $\Delta\tilde{\theta} = \Delta\theta$; if not, without loss of generality, assume that $A$ is in the first quadrant with $\pi/4 \leq \theta_A < \pi/2$, then $\tan\theta_{\tilde{A}} = \frac{b}{a\mu_1^2} = \frac{\tan\theta_A}{\mu_1^2}$, and by (2.18), we have

$$\Delta\tilde{\theta} \geq \theta_{\tilde{A}} \geq \frac{\theta_A}{4\mu_1^2} \geq \frac{\Delta\theta}{16\mu_1^2}$$

($\theta_A \geq \Delta\theta/4$ since $\theta_A \geq \pi/4$) and then $\frac{\Delta\theta}{16\mu_1^2} \leq \Delta\tilde{\theta} \leq 2\Delta\theta$.

*Case 2.* $0 < \Delta\theta < \pi/2$. In this case, $(I, III)$ and $(II, IV)$ are impossible. We will divide the following proof into three parts.

(1). We discuss $(I, I)$ here; the argument may be readily adapted to $(II, II)$, $(III, III)$, and $(IV, IV)$. Without loss of generality, assume $\Delta\theta = \theta_A - \theta_B$, then we get

$$\tan\Delta\tilde{\theta} = \frac{\mu_1^2(bc - ad)}{bd + \mu_1^4 ac} \geq \frac{\tan\Delta\theta}{\mu_1^2},$$

and by (2.18), we get $\frac{\Delta\theta}{4\mu_1^2} \leq \Delta\tilde{\theta}$. Similarly, we will get $\Delta\tilde{\theta} \leq 4\mu_1^2\Delta\theta$ since $\tan\Delta\tilde{\theta} \leq \mu_1^2\tan\Delta\theta$, and so $\frac{\Delta\theta}{4\mu_1^2} \leq \Delta\tilde{\theta} \leq 4\mu_1^2\Delta\theta$ follows.

(2). We discuss $(I, IV)$ here; an adaptation handles $(II, III)$. Without loss of generality, assume $\theta_A \geq \Delta\theta/2$. Obviously, we have $\Delta\tilde{\theta} \leq \Delta\theta$. Conversely, we have $\frac{\Delta\theta}{16\mu_1^2} \leq \Delta\tilde{\theta}$ (it is essentially the same as $(I, IV)$ in *Case 1*), and so $\frac{\Delta\theta}{16\mu_1^2} \leq \Delta\tilde{\theta} \leq \Delta\theta$ follows.
(3). We discuss $(I, II)$ here, and the method can be applied to $(III, IV)$. Obviously, we have $\Delta\theta \leq \Delta\tilde{\theta}$. Without loss of generality, assume that $A$ is in the first quadrant and makes an angle $h_A$ with the $y$-axis and that $B$ is in the second quadrant and makes an angle $h_B$ with the $y$-axis. Clearly, $\Delta\theta = h_A + h_B$. Let $h_{\tilde{A}}$ and $h_{\tilde{B}}$ be the angle between the $y$-axis and $\tilde{A}$ and the angle between the $y$-axis and $\tilde{B}$, respectively. By (2.18), we conclude that $h_{\tilde{A}} \leq 4\mu_1^2 h_A$ since $\tan h_{\tilde{A}} = \mu_1^2 \tan h_A$. Similarly, we get $h_{\tilde{B}} \leq 4\mu_1^2 h_B$. So it follows that $\Delta\theta \leq \Delta\tilde{\theta} = h_{\tilde{A}} + h_{\tilde{B}} \leq 4\mu_1^2(h_A + h_B) = 4\mu_1^2\Delta\theta$.

Through the above analysis, we see that $\frac{\Delta\theta}{16\mu_1^2} \leq \Delta\tilde{\theta} \leq 16\mu_1^2\Delta\theta$, concluding the proof of Claim 2.31. □

By this claim, we can modify the last paragraph of the proof of [1, Claim 3.8] as stated below and then our lemma follows.

"If $r$ sufficiently large, we conclude that for every $0 \leq l \leq r-2$, there exists a rotation $R_{l,j}$ of angle $\theta_j$ with $r^{-2.5N} < \theta_j < r^{-0.3N}$ such that $C^{(l+1,i_j)}u_j = R_{l,j}C^{(l,i_j)}u_j$. It immediately follows that there exists at most one choice of $0 \leq t_{i_j} \leq r-1$ such that $C^{(t_{i_j},i_j)}u_j$ has angle at most $r^{-3N}$ with $s_{j+1}$, as desired."

We would like to explain how to obtain the statement described in the paragraph above. If $r$ is sufficiently large, it is not hard to conclude that for every $0 \leq l \leq r-2$, there exists a rotation $\tilde{R}_{l,j}$ of angle $\tilde{\theta}_j$ with $r^{-2N} < \tilde{\theta}_j < r^{-0.5N}$ such that $A_{n_k}^{(E-\lambda r^{-N}(l+1),\lambda f_{i_j})}u_j = \tilde{R}_{l,j}A_{n_k}^{(E-\lambda r^{-N}l,\lambda f_{i_j})}u_j$. If $i_d = m$, we have

$$\begin{aligned} C^{(l+1,m)}u_m &= A_{n_k}^{(E,\lambda f_m)}A_{n_k}^{(E-\lambda r^{-N}(l+1),\lambda f_m)}u_m \\ &= A_{n_k}^{(E,\lambda f_m)}\tilde{R}_{l,m}A_{n_k}^{(E-\lambda r^{-N}(l),\lambda f_m)}u_m. \end{aligned} \tag{2.19}$$

Since $A_{n_k}^{(E,\lambda f_m)} \in \mathrm{SL}(2,R)$ is independent of $r$, we can apply Claim 2.31 to (2.19) so that we have

$$\begin{aligned} C^{(l+1,m)}u_m &= A_{n_k}^{(E,\lambda f_m)}\tilde{R}_{l,m}A_{n_k}^{(E-\lambda r^{-N}(l),\lambda f_m)}u_m \\ &= R_{l,m}A_{n_k}^{(E,\lambda f_m)}A_{n_k}^{(E-\lambda r^{-N}(l),\lambda f_m)}u_m \\ &= R_{l,m}C^{(l,m)}u_m, \end{aligned}$$

where $R_{l,m}$ is a rotation of angle $\theta_m$ with $r^{-2.5N} < \theta_m < r^{-0.3N}$. Then the above paragraph follows. $\square$

Now we can give the

*Proof of Theorem 2.23.* Given a $p_0$-periodic $f \in C(\Omega, \mathbb{R})$ and $0 < \varepsilon_0 < 1$, consider $B_{\varepsilon_0}(f) \subset C(\Omega, \mathbb{R})$. (We will work within this ball. The denseness of periodic potentials then implies the denseness of our constructed limit-periodic potentials.) Let $N$ from Lemma 2.30 be 2. Let $\varepsilon_1 = \frac{\varepsilon_0}{10}$. By Lemma 2.29, there exists a finite family $F_1 = \{f_1, f_2, \ldots, f_{m_1}\}$ of $p_1$-periodic sampling functions such that $F_1 \subset B_{\varepsilon_0}(f)$ and $L(E, \lambda F_1) > \delta_1$ for some $0 < \delta_1 < 1$ whenever $\varepsilon_1 < |\lambda| < \frac{1}{\varepsilon_1}$ and $E \in \mathbb{R}$ (note that $L(E, \lambda f_i) \geq 1$ if $|E| \geq \|\lambda f_i\| + 4$). Our constructions start with $F_1$ and we will divide them into several steps.

*Construction 1.* First, we will apply Lemma 2.29 to $F_1$ in order to enlarge the range of $\lambda$'s. Let $\varepsilon_2 = \frac{min\{\varepsilon_1, \delta_1\}}{10}$. Then, there exists a finite family of $\tilde{p}_1$-periodic potentials $\tilde{F}_1 = \{\tilde{f}_1, \tilde{f}_2, \ldots, \tilde{f}_{\tilde{m}_1}\} \subset B_{\varepsilon_0}(f)$ such that

$$L(E, \lambda \tilde{F}_1) > \tilde{\delta}_1$$

for some $0 < \tilde{\delta}_1 < 1$ whenever $\forall \varepsilon_2 < |\lambda| < \frac{1}{\varepsilon_2}$ and $E \in \mathbb{R}$, and

$$\left| L(E, \lambda \tilde{F}_1) - L(E, \lambda F_1) \right| < \frac{\varepsilon_2}{2} \tag{2.20}$$

whenever $|E| < \frac{1}{\varepsilon_2}$ and $|\lambda| < \frac{1}{\varepsilon_2}$.

Explicitly, the construction of $\tilde{F}_1$ follows from the proof of [1, Claim 3.1]. For very large $\tilde{p}_1 > p_1$ (it must obey $p_1 | \tilde{p}_1$), choose $N_1(\tilde{p}_1)$ such that if $|E| < \frac{1}{\varepsilon_2}$, $|\lambda| \leq \frac{1}{\varepsilon_2}$, $f_i \in F_1$ and a $\tilde{p}_1$-periodic potential $\tilde{f}$ which is $\frac{2p_1+1}{N_1(\tilde{p}_1)}$ close to $f$ then $|L(E, \lambda \tilde{f}) - L(E, \lambda f)| < \frac{\varepsilon_2}{2}$, since the Lyapunov exponent is continuous for periodic potentials (see Proposition 2.11).

For $1 \leq j \leq 2p_1 + 1$, we define $\tilde{p}_1$-periodic potentials $\tilde{f}^{(i,j)}$ by $\tilde{f}^{(i,j)}(n) = f_i(n), 0 \leq n \leq \tilde{p}_1 - 2$ and $\tilde{f}^{(i,j)}(\tilde{p}_1 - 1) = f_i(\tilde{p}_1 - 1) + \frac{j}{N_1(\tilde{p}_1)}$. By [1, Claim 3.4], there exists $j_0$ such that the spectrum of $\tilde{f}^{(i,j_0)}$ has exactly $\tilde{p}_1$ components, that is, all gaps of its spectrum are open. For convenience, we write $\tilde{f}^{(i)} = \tilde{f}^{(i,j_0)}$. So there exists $h = h(F_1, \tilde{p}_1, \varepsilon_2) > 0$ such that for any $f_i \in F_1$ and $\varepsilon_2 \leq |\lambda| \leq \frac{1}{\varepsilon_2}$, $\Sigma(\lambda \tilde{f}^{(i)})$ has $\tilde{p}_1$ components and the Lebesgue measure of the smallest gap is at least $h$. Choose an integer $N_2(\tilde{p}_1)$ with $N_2(\tilde{p}_1) > \frac{4\pi}{\varepsilon_2 h \tilde{p}_1}$.

For $0 \leq l \leq N_2(\tilde{p}_1)$, let $\tilde{f}^{(i,l)} = \tilde{f}^{(i)} + \frac{4\pi l}{\varepsilon_2 \tilde{p}_1 N_2(\tilde{p}_1)}$. Then $\tilde{F}_1$ is just the family obtained by collecting the $\tilde{f}^{(i,l)}$ for different $f_i \in F_1$ and $0 \leq l \leq N_2(\tilde{p}_1)$. Order $\tilde{F}_1$ as $\tilde{F}_1 = \{\tilde{f}_1, \tilde{f}_2, \ldots, \tilde{f}_{\tilde{m}_1}\}$ such that $\tilde{f}_1 = \tilde{f}^{(1,0)}$ and $\tilde{f}_{\tilde{m}_1} = \tilde{f}^{(1,1)}$. We can also assume that $N_2(\tilde{p}_1)$ was chosen large enough, so that we have $\|\tilde{f}_{\tilde{m}_1} - \tilde{f}_1\| = \frac{4\pi}{\varepsilon_2 \tilde{p}_1 N_2(\tilde{p}_1)} < 1/3$ (this will be used to conclude that our limit-periodic potentials are Gordon potentials).

*Construction 2.* Applying Lemma 2.30 to $\tilde{F}_1$, there exists a finite family of $p_2$-periodic potentials $F_2 = \{f_2^{\vec{t}_1}, f_2^{\vec{t}_2}, \ldots, f_2^{\vec{t}_{m_2}}\}$ such that

$$F_2 \subset B_{p_2^{-2}} \subset B_{\varepsilon_2} \subset B_{\varepsilon_0}(f)$$

and

$$L(E, \lambda F_2) > \delta_2$$

for some $0 < \delta_2 < 1$ whenever $\varepsilon_2 < |\lambda| < \frac{1}{\varepsilon_2}$ and $E \in \mathbb{R}$, and

$$\left| L(E, \lambda F_2) - L(E, \lambda \tilde{F}_1) \right| < \frac{\varepsilon_2}{2} \tag{2.21}$$

whenever $|E|$, $|\lambda| < \frac{1}{\varepsilon_2}$. From (2.20) and (2.21), we have

$$|L(E, \lambda F_2) - L(E, \lambda F_1)| < \varepsilon_2$$

for $|E|$, $|\lambda| < \frac{1}{\varepsilon_2}$.

Explicitly, we construct $F_2$ as follows (cf. the proof of Lemma 2.30). Let $p_2$ large and $p_2 = \tilde{m}_1 \tilde{p}_1 r_2 + d$, $0 \leq d \leq \tilde{m}_1 \tilde{p}_1 - 1$. Let $I_j = [j\tilde{p}_1, (j+1)\tilde{p}_1 - 1] \subset \mathbb{Z}$ and let $0 = j_0 < j_1 < \cdots < j_{\tilde{m}_1 - 1} < j_{\tilde{m}_1} = \frac{p_2}{\tilde{p}_1}$ be a sequence such that $j_{i+1} - j_i = r_2 + 1$ when $0 \leq i < d$ and $j_{i+1} - j_i = r_2$ when $d \leq i \leq \tilde{m}_1 \tilde{p}_1 - 1$. Define a $p_2$-periodic potential $f_2(l)$ for $0 \leq l \leq p_2 - 1$ as follows. Let $j$ be such that $l \in I_j$ and let $i$ be such that $j_{i-1} \leq j < j_i$ and let $f_2(l) = \tilde{f}_i(l)$. For any sequence $\vec{t} = (t_1, t_2, \ldots, t_{\tilde{m}_1})$ with $t_i \in \{0, 1, \ldots, r_2 - 1\}$, let $f_2^{\vec{t}}$ be a $p_2$-periodic potential defined as follows. Let $0 \leq l \leq p_2 - 1$, and let $j$ be such that $l \in I_j$. If $j = j_i - 1$ for some $1 \leq i < \tilde{m}_1$, we let $f_2^{\vec{t}}(l) = f_2(l) + r_2^{-4} t_i$, and $j = j_{\tilde{m}_1} - 2$ then $f_2^{\vec{t}}(l) = f_2(l) + r_2^{-4} t_{\tilde{m}_1}$. Otherwise we let $f_2^{\vec{t}}(l) = f_2(l)$. Let $p_2$ be sufficiently large so that $p_2^{-2} < 1/3$.

Moreover, we can estimate the Lebesgue measure of the spectrum. For any $E \in \mathbb{R}$ and $\varepsilon_2 < |\lambda| < \frac{1}{\varepsilon_2}$, we can find $\tilde{f}_i \in \tilde{F}_1$ such that $L(E, \lambda \tilde{f}_i) > \tilde{\delta}_1$ since $L(E, \lambda \tilde{F}_1) > \tilde{\delta}_1$. If $r_2$ large enough, we have $\|A_{(r_2-2)\tilde{p}_1}^{(E,\lambda \tilde{f}_i)}\| > e^{\tilde{\delta}_1 (r_2 - 2)\tilde{p}_1}$. Then we have

$$\|A_{(r_2-2)\tilde{p}_1}^{(E,\lambda f_2^{\vec{t}_k})}(f_2^{\vec{t}_k}(j_{i-1}\tilde{p}_1))\| = \|A_{(r_2-2)\tilde{p}_1}^{(E,\lambda \tilde{f}_i)}\| > e^{\tilde{\delta}_1 (r_2-2)\tilde{p}_1}.$$

Since $E$ is arbitrary, we can apply Lemma 2.28 to conclude that the total Lebesgue measure of $\Sigma(\lambda f_2^{\vec{t}_k})$ is at most $4\pi p_2 e^{-\tilde{\delta}_1(r_2-2)\tilde{p}_1} < e^{-\tilde{p}_1 p_2^{1/2}}$ when $r_2$ sufficiently large. (Here $f_2^{\vec{t}_k}$ can be any element from $F_2$.)

*Construction 3.* Repeating the above procedures. Once we have constructed $F_{i-1} \subset B_{p_{i-1}^{-(i-1)}} \subset B_{\varepsilon_{i-1}}$, by Lemma 2.29, we can get a finite family of $\tilde{p}_{i-1}$-periodic potentials $\tilde{F}_{i-1} \subset B_{p_{i-1}^{-(i-1)}}$ satisfying the following. Let $\varepsilon_i = \frac{min\{\varepsilon_{i-1},\delta_{i-1}\}}{10}$, and we have

$$L(E, \lambda\tilde{F}_{i-1}) > \tilde{\delta}_{i-1}$$

for some $0 < \tilde{\delta}_{i-1} < 1$ whenever $\forall\varepsilon_i < |\lambda| < \frac{1}{\varepsilon_i}$ and $E \in \mathbb{R}$, and

$$\Big|L(E, \lambda\tilde{F}_{i-1}) - L(E, \lambda F_{i-1})\Big| < \frac{\varepsilon_i}{2}$$

whenever $|E| < \frac{1}{\varepsilon_i}$ and $|\lambda| < \frac{1}{\varepsilon_i}$.

Next, as in *Construction 2*, we will get a finite family $F_i$ of $p_i$-periodic potentials which satisfies the following (here our perturbation is $r_i^{-Ni}t = r_i^{-2i}t$, $t \in \{0, 1, 2, \ldots, r_i - 1\}$).

(i). $L(E, \lambda F_i) > \delta_i$ for some $0 < \delta_i < 1$ and all $E \in \mathbb{R}$ and $\varepsilon_i < |\lambda| < \varepsilon_i^{-1}$.

(ii). $|L(E, \lambda F_i) - L(E, \lambda F_{i-1})| < \varepsilon_i$, for $|E| < \frac{1}{\varepsilon_i}$ and $|\lambda| < \frac{1}{\varepsilon_i}$.

(iii). $F_i \subset B_{p_i^{-i}} \subset B_{\varepsilon_i} \subset B_{\varepsilon_{i-1}} \subset B_{\varepsilon_0}(f), i > 2$. (Note: $B_{\varepsilon_2}$ may not be in $B_{\varepsilon_1}$.)

(iv). $\forall f_i^{\vec{t}_k} \in F_i$, $|\Sigma(\lambda f_i^{\vec{t}_k})| \leq e^{-\tilde{p}_{i-1}p_i^{1/2}}$ (here $|\cdot|$ denotes the Lebesgue measure) when $\varepsilon_i < |\lambda| < \varepsilon_i^{-1}$.

(v). $p_i^{-i} < \frac{1}{3}(i-1)^{-\tilde{p}_{i-1}}$ since we can let $p_i$ be sufficiently large.

(vi). $\|f_i^{\vec{t}_1} - f_i^{\vec{t}_{m_i}}\| = \frac{4\pi}{\varepsilon_i\tilde{p}_{i-1}N_2(\tilde{p}_{i-1})} < \frac{1}{3}(i-1)^{-\tilde{p}_{i-1}}$. Here $N_2(\tilde{p}_{i-1})$ appears as in *Construction 1*, and we can ensure that this inequality holds since $\tilde{p}_{i-1}$ is fixed while $N_2(\tilde{p}_{i-1})$ can be taken as large as needed.

Then we will get a limit-periodic potential $f_\infty \in B_{\varepsilon_0}(f)$, whose Lyapunov exponent is a positive continuous function of energy $E$ and the Lebesgue measure of the spectrum is zero (Lemma 2.12 implies that $L(E, \lambda f_i^{\vec{t}}) \to L(E, \lambda f_\infty)$). Moreover, we have the following two claims.

**Claim 2.32.** *$f_\infty$ is a Gordon potential.*

*Proof.* Let $q_i = \tilde{p}_i$. Obviously, $q_i \to \infty$ as $i \to \infty$. For $i \geq 1$, we have

$$\begin{aligned}
\max_{1\leq n\leq q_i} |f_\infty(n) - f_\infty(n \pm q_i)| &\leq |f_\infty(n) - f_{i+1}^{\vec{t}_1}(n)| + |f_\infty(n \pm q_i) - f_{i+1}^{\vec{t}_1}(n \pm q_i)| \\
&\quad + |f_{i+1}^{\vec{t}_1}(n) - f_{i+1}^{\vec{t}_1}(n \pm q_i)| \\
&\leq p_{i+1}^{-(i+1)} + p_{i+1}^{-(i+1)} + \frac{4\pi}{\varepsilon_{i+1}\tilde{p}_i N_2(\tilde{p}_i)} \\
&\leq 2\frac{1}{3}(i)^{-\tilde{p}_i} + \frac{1}{3}(i)^{-\tilde{p}_i} \\
&\leq i^{-q_i}.
\end{aligned}$$

So $f_\infty$ is a Gordon potential. (Here $f_{i+1}^{\vec{t}_1}$ is an element of $F_{i+1}$). □

**Claim 2.33.** *$\Sigma(\lambda f_\infty)$ has zero Hausdorff dimension for every $\lambda \neq 0$.*

*Proof.* Let $\lambda \neq 0$ and $0 < \alpha \leq 1$ be given. Without loss of generality, assume $\lambda > 0$. Choose $i$ large enough so that $\varepsilon_i < \lambda < 1/\varepsilon_i$ and $1/i < \alpha$. For every $f_i^{\vec{t}_k} \in F_i$, $\|\lambda f_\infty - \lambda f_i^{\vec{t}_k}\| < \lambda p_i^{-i}$ implies[2] $\operatorname{dist}(\Sigma(\lambda f_\infty), \Sigma(\lambda f_i^{\vec{t}_k})) < \lambda p_i^{-i}$. Since $\lambda f_i^{\vec{t}_k}$ is $p_i$-periodic, we have

$$\Sigma(\lambda f_i^{\vec{t}_k}) = \bigcup_{z=1}^{p_i} \tilde{I}_z^{(\vec{t}_k,i)},$$

where $\tilde{I}_z^{(\vec{t}_k,i)} = [a_z, b_z]$ is a closed interval.

Let $I_z^{(\vec{t}_k,i)} = [a_z - \lambda p_i^{-i}, b_z + \lambda p_i^{-i}]$ and since $\operatorname{dist}(\Sigma(\lambda f_\infty), \Sigma(\lambda f_i^{\vec{t}_k})) \leq \lambda p_i^{-i}$, we have

$$\Sigma(\lambda f_\infty) \subset \bigcup_{z=1}^{p_i} I_z^{(\vec{t}_k,i)}.$$

Moreover, $b_z - a_z \leq e^{-\tilde{p}_{i-1} p_i^{1/2}}$ since $|\Sigma(\lambda f_i^{\vec{t}_k})| \leq e^{-\tilde{p}_{i-1} p_i^{1/2}}$. Then we have

$$\begin{aligned}
h^\alpha(\Sigma(\lambda f_\infty)) &\leq \lim_{i\to\infty} \sum_{z=1}^{p_i} (e^{-\tilde{p}_{i-1} p_i^{1/2}} + 2\lambda p_i^{-i})^\alpha \\
&= \lim_{i\to\infty} p_i (e^{-\tilde{p}_{i-1} p_i^{1/2}} + 2\lambda p_i^{-i})^\alpha \\
&= \lim_{i\to\infty} (p_i^{1/\alpha} e^{-\tilde{p}_{i-1} p_i^{1/2}} + 2\lambda p_i^{-i+1/\alpha})^\alpha.
\end{aligned}$$

[2]It is well known that for $V, W : \mathbb{Z} \to \mathbb{R}$ bounded, we have $\operatorname{dist}(\sigma(\Delta + V), \sigma(\Delta + W)) \leq \|V - W\|_\infty$, where $\operatorname{dist}(A, B)$ denotes the Hausdorff distance of two compact subsets $A, B$ of $\mathbb{R}$.

Since $1/i < \alpha$, we have $-i + 1/\alpha < 0$, and it follows that

$$\lim_{i\to\infty} (p_i^{1/\alpha} e^{-p_{i-1}p_i^{1/2}} + 2\lambda p_i^{-i+1/\alpha})^\alpha = 0.$$

So we have $h^\alpha(\Sigma(\lambda f_\infty)) = 0$ (note: when $i \to \infty$, $\lambda$ belongs to $(\varepsilon_i, \frac{1}{\varepsilon_i})$ for all $i$ large enough since this interval is expanding). So the Hausdorff dimension of the spectrum is less than $\alpha$. Since $\alpha$ was arbitrary, the Hausdorff dimension must be zero. □

This implies all the assertions in Theorem 2.23 except for the absence of eigenvalues for every $\omega$. Given the Gordon Lemma (see Lemma 2.20 above), this last statement will follow once Theorem 2.25 is established. □

**Remark 2.34.** *Since $\delta_i \le \delta_{i-1}/10, i \ge 1$, it is true that when $\varepsilon_i < |\lambda| < \frac{1}{\varepsilon_i}$, $L(E, \lambda f_\infty) \ge \frac{8}{9}\delta_i$ for any $E \in \mathbb{R}$. This gives information about the range of the Lyapunov exponent on certain intervals. Clearly, $\delta_i \to 0$ when $i \to \infty$ since the Lyapunov exponent will go to zero when $\lambda$ goes to zero.*

*Proof of Theorem 2.25.* Let $\omega = e$ first. Relative to any minimal translation $\tilde{T}$, the selected $f$ in the proof of Theorem 2.23 is still $n_0$-periodic by Proposition 2.10, so we can start with the same ball $B_{\varepsilon_0}(f)$ and choose the same periodic potentials in $B_{\varepsilon_0}(f)$. Then we get the same $f_\infty$. For the finite family $F_i$ from *Construction 3*, though the Lyapunov exponent may change, the following properties hold (note that $\|f_i^{\vec{t}_1}\|$ does not change).

(i). $F_i \subset B_{p_i^{-i}} \subset B_{\varepsilon_i} \subset B_{\varepsilon_0}(f)$.

(ii). $p_i^{-i} < \frac{1}{3}(i-1)^{-\tilde{p}_{i-1}}$.

(iii). $\|f_i^{\vec{t}_1} - f_i^{\vec{t}_{m_i}}\| = \frac{4\pi}{\varepsilon_i \tilde{p}_{i-1} N_2(\tilde{p}_{i-1})} < \frac{1}{3}(i-1)^{-\tilde{p}_{i-1}}$.

Then Claim 2.32 holds true, and so $f(\tilde{T}^n(e))$ is a Gordon potential. For arbitrary $\tilde{\omega}$, if we repeat the same procedures, (i)—(iii) above still hold as stated (since none of them are related to $\omega$), and Theorem 2.25 follows. □

## 2.4 Uniform Localization

Localization is a topic that has been explored in the context of Schrödinger operators to a great extent. By now several mechanisms are known that lead to localization, at least in suitable energy regions. The most important one is randomness or, more generally, weak correlations. This aspect goes back to the seminal paper [2] of Anderson. Another important mechanism is strong coupling and, related to this, positive Lyapunov exponents. The latter approach can be used to prove localization for strongly correlated potentials.

While localization does not occur for periodic potentials. Limit-periodic potentials are closest to periodic potentials (at least among the stationary ones) and hence for them, one would expect either the absence of localization or a difficult localization proof in the rare cases where it holds. There are two notable examples about pure point spectrum in the study of limit-periodic Schrödinger operators. The first is a paper by Chulaevsky and Molchanov, [29], which unfortunately does not contain a proof of the theorem on the presence of pure point spectrum for some continuum one-dimensional limit-periodic Schrödinger operators stated there. Moreover, their examples have zero Lyapunov exponent and hence are not localized in the standard sense. The other relevant paper is Pöschel's work [30], where he proves a general theorem that provides a sufficient condition for uniform localization along with two examples showing that the general result is applicable to limit-periodic potentials. Incidentally, the Chulaevsky-Molchanov paper uses features of randomness while Pöschel's paper uses strong coupling.

In this section, we will prove uniform localization results that hold uniformly for all elements of the hull $\Omega$. This is a novel phenomenon. Indeed, usually localization can be proved, and in fact holds, only almost surely. For random potentials, this is obvious since there are periodic realizations of the potential. For certain almost periodic potentials, there are results to this effect due to Jitomirskaya-Simon [26] and Gordon [22]. Regarding results establishing pure point spectrum for all elements of the family, we are aware of the following: For the Maryland model, see [20, 19, 23, 32, 39], which has an unbounded potential (and hence is not almost periodic), pure point spectrum was shown for the whole family but without uniform decay of eigenfunctions. There is some unpublished work of Jitomirskaya establishing a similar result for a bounded non-almost periodic model. To the best of our knowledge, we exhibit the first almost periodic example that is uniformly localized across the hull and the spectrum. The result has been published in Journal d'Analyse Mathematique [14].

### 2.4.1 Introduction to Uniform Localization

**Definition 2.35.** We say that a family $\{u_k\} \subset \ell^2(\mathbb{Z})$ is uniformly localized if there exist constants $r > 0$, called the decay rate, and $c < \infty$ such that for every element $u_k$ of the family, one can find $m_k \in \mathbb{Z}$, called the center of localization, so that $|u_k(n)| \leq ce^{-r|n-m_k|}$ for every $n \in \mathbb{Z}$. We say that the operator $H_\omega$ has ULE if it has a complete set of uniformly localized eigenfunctions.[3]

The notion of uniformly localized eigenfunctions and related ones were introduced by del Rio et al. in their comprehensive study of the question "What is localization?" [15, 16]. As explained there, ULE implies uniform dynamical localization, that is, if $H_\omega$ has ULE, then

$$\sup_{t\in\mathbb{R}} \left|\left\langle \delta_n, e^{-itH_\omega}\delta_m \right\rangle\right| \leq C_\omega e^{-r_\omega|n-m|} \tag{2.22}$$

with suitable constants $C_\omega, r_\omega \in (0, \infty)$. While both properties are desirable, they are extremely rare. To quote from [16], "the problem is that ULE does not occur" and "it is an open question, in fact, whether there is *any* Schrödinger operator with ULE." Del Rio et al. may not have been aware of Pöschel's work [30] since it predates theirs and provides some examples of Schrödinger operators with ULE.

The occurrence of pure point spectrum for the operators $\{H_\omega\}_{\omega\in\Omega}$ is called *phase stable* if it holds for every $\omega \in \Omega$. It is an unusual phenomenon since most known models are not phase stable. It is known that uniform localization of eigenfunctions (ULE) has a close connection with phase stability of pure point spectrum; compare the following theorem.

**Theorem 2.36.** [15, Theorem C.1] *If $H_\omega$ has ULE for $\omega$ in a set of positive $\mu$-measure, then $H_\omega$ has pure point spectrum for every $\omega \in \operatorname{supp}(\mu)$, where* $\operatorname{supp}(\mu)$ *is the complement of the largest open set $S \subset \Omega$ for which $\mu(S) = 0$.*

Jitomirskaya pointed out in [25] that Theorem 2.36 can be strengthened for a minimal $T$ in the sense that if there exists some $\omega_0$ such that $H_{\omega_0}$ has ULE, then $H_\omega$ has pure point spectrum for every $\omega \in \operatorname{supp}(\mu)$.

The hull of a limit-periodic potential is a procyclic group, so that we can classify such procyclic groups by their frequency integer sets. Every procyclic group has a unique maximal frequency integer set $S = \{n_k\} \subseteq \mathbb{Z}_+$ with the property that $n_{k+1}/n_k$ is prime for every $k$.

[3] Recall that a set of vectors is called complete if their span (i.e., the set of finite linear combinations of vectors from this set) is dense.

**Definition 2.37.** For a procyclic group, we say that it satisfies the condition $\mathscr{A}$ if its maximal frequency integer set $S = \{n_k\} \subseteq \mathbb{Z}_+$ has the following property: there exists some integer $m \geq 2$ such that for every $k$, we have $n_k < n_{k+1} \leq n_k^m$, that is, $\log n_{k+1}/\log n_k$ is uniformly bounded.

We can now state our main result, which is just [14, Theorem 2.5.].

**Theorem 2.38** (Uniform Localization)**.** *Suppose $\Omega$ is a procyclic group with the condition $\mathscr{A}$, and $T : \Omega \to \Omega$ a minimal translation. Then there exists some $f \in C(\Omega, \mathbb{R})$ such that for every $\omega \in \Omega$, the Schrödinger operator with potential $f(T^n(\omega))$ has ULE with $\omega$-independent constants. In particular, we have uniform dynamical localization* (2.22) *for every $\omega$ with $\omega$-independent constants as well.*

We will heavily use Pöschel's results in [30], which will be recalled in Subsection 2.4.3, to obtain the above theorem. Pöschel used an abstraction of KAM methods, with some of the basic ideas going back to Craig [9], Rüssmann [35] and Moser [28]. In this approach, there is an important concept, that of a distal sequence, which we will discuss in Subsection 2.4.2. The first step in proving Theorem 2.38 is to construct a distal limit-periodic potential in our framework.

### 2.4.2 Distal Sequences

Here we will discuss approximation functions and distal sequences; compare [30] and [35].

**Definition 2.39.** A function $Q(x) : [0, \infty) \to [1, \infty)$ is called an approximation function if both

$$q(t) = t^{-4} \sup_{x \geq 0} Q(x) e^{-tx}$$

and

$$h(t) = \inf_{\kappa_t} \prod_{i=0}^{\infty} q(t_i)^{2^{-i-1}} \tag{2.23}$$

are finite for every $t > 0$. In (2.23), $\kappa_t$ denotes the set of all sequences $t \geq t_1 \geq t_2 \geq \cdots \geq 0$ with $\sum t_i \leq t$.

**Definition 2.40.** A sequence $V \in \ell^\infty(\mathbb{Z})$ is called distal if for some approximation function $Q$, we have

$$\inf_{i\in\mathbb{Z}} |V_i - V_{i+k}| \geq Q(|k|)^{-1}$$

for every $k \in \mathbb{Z} \setminus \{0\}$.

**Proposition 2.41.** *If $V \in \ell^\infty(\mathbb{Z})$ is distal, then every $\tilde{V} \in \mathrm{hull}(V)$ is also distal.*

*Proof.* This follows readily from the definition. □

The following lemma shows how to generate distal sequences in our framework.

**Lemma 2.42.** *Suppose $\Omega$ is a procyclic group with the condition $\mathscr{A}$, and $T : \Omega \to \Omega$ a minimal translation. There exists an $f \in C(\Omega, \mathbb{R})$ such that $(f(T^i(e)))_{i\in\mathbb{Z}}$ is a distal sequence.*

*Proof.* Given a procyclic group $\Omega$ and a minimal translation $T$, by Lemma 1.33 there is a limit-periodic potential $L$ such that $\mathrm{hull}(L) \cong \Omega$. Since $\Omega$ satisfies the condition $\mathscr{A}$, there exists $m \geq 2$ such that for the elements of its maximal frequency integer set $S_\Omega = \{n_k\}$, we have $n_{k-1} < n_k \leq n_{k-1}^m$ for every $k$.

Consider $S_\Omega$. Here we let $n_1 > 1$. For $n_1 \in S_\Omega$, there must exist some $n_k \in [n_1^3, n_1^{3m}]$. If not, we pick the largest $n_i \in [n_1, n_1^3)$ and then $n_{i+1}$ will be strictly larger than $n_1^{3m}$. Then we have $n_{i+1} > n_1^{3m} > n_i^m$ which contradicts the assumption. So we can pick $n_k$ such that $n_1^3 \leq n_k \leq n_1^{3m}$. By induction, we can pick a subset of $S_\Omega$ which we still denote by $I_0 = \{n_k\}$ satisfying $n_k^3 \leq n_{k+1} \leq n_k^{3m}$ for every $k \in \mathbb{Z}^+$. Without any contradiction, we take $n_0 = 1$ for the following computation.

Define $a_v(i) = j$ where $0 \leq j < n_v$ and $i = j \pmod{n_v}$, so $a_v$ is $n_v$-periodic. Let $V = (V_i)_{i\in\mathbb{Z}}$ and $V^{(k)} = (V_i^{(k)})_{i\in\mathbb{Z}}$, where

$$V_i = \sum_{v=1}^{\infty} \frac{a_v(i)}{n_{v-1}^2 n_v} \quad \text{and} \quad V_i^{(k)} = \sum_{v=1}^{k} \frac{a_v(i)}{n_{v-1}^2 n_v}.$$

By the divisibility property of any frequency integer set, $V^{(k)}$ is an $n_k$-periodic sequence. Since for every $i \in \mathbb{Z}$ and $k \in \mathbb{Z}_+$ we have

$$\left|V_i - V_i^{(k)}\right| = \left|\sum_{v=k+1}^{\infty} \frac{a_v(i)}{n_{v-1}^2 n_v}\right| \leq \sum_{v=k+1}^{\infty} \frac{1}{n_{v-1}^2},$$

it follows that $V^{(k)}$ converges to $V$ uniformly. Thus, $V$ is limit-periodic and one of its frequency integer sets is $I_0$.

For any $i_1 \neq i_2$, fix $k$ so that $n_{k-1} \leq |i_1 - i_2| < n_k$. If $k = 1$, then $|V_{i_1}^{(1)} - V_{i_2}^{(1)}| \geq \frac{1}{n_1}$. Also, we have

$$\begin{aligned} &|(V_{i_1} - V_{i_1}^{(1)}) - (V_{i_2} - V_{i_2}^{(1)})| \\ &\leq n_1 \sum_{v=2}^{\infty} \frac{1}{n_{v-1}^2 n_v} \\ &\leq \frac{8}{7n_1 n_2} \\ &\leq \frac{4}{7n_1} \end{aligned}$$

So it is easy to see that $|V_{i_1} - V_{i_2}| \geq \frac{3}{7n_1} \geq \frac{2}{3n_1^{3m+1}}$.

If $k \geq 2$, we have

$$\frac{1}{n_{k-2}^2 n_{k-1}} > \frac{2(n_k - 1)}{n_{k-1}^2 n_k}, \tag{2.24}$$

since $n_i^3 \leq n_{i+1} \leq n_i^{3m}$. Moreover, we have

$$\begin{aligned} |V_{i_1}^{(k)} - V_{i_2}^{(k)}| &= \left| \sum_{v=1}^{k} \frac{(a_v(i_1) - a_v(i_2))}{n_{v-1}^2 n_v} \right| \\ &= \left| \sum_{v=1}^{k-1} \frac{(a_v(i_1) - a_v(i_2))}{n_{v-1}^2 n_v} + \frac{i_1 - i_2}{n_{k-1}^2 n_k} \right|. \end{aligned}$$

From (2.24) we conclude that $|\sum_{v=1}^{k-1} \frac{(a_v(i_1)-a_v(i_2))}{n_{v-1}^2 n_v}|$ is 0 or larger than $\frac{2(n_k-1)}{n_{k-1}^2 n_k}$. So $|V_{i_1}^{(k)} - V_{i_2}^{(k)}| = \left|\sum_{v=1}^{k-1} \frac{(a_v(i_1)-a_v(i_2))}{n_{v-1}^2 n_v} + \frac{i_1-i_2}{n_{k-1}^2 n_k}\right| \geq \frac{n_{k-1}}{n_{k-1}^2 n_k} = \frac{1}{n_{k-1} n_k}$. We also

have

$$\begin{aligned}
&|(V_{i_1} - V_{i_1}^{(k)}) - (V_{i_2} - V_{i_2}^{(k)})| \\
&= \left| \sum_{v=k+1}^{\infty} \frac{(a_v(i_1) - a_v(i_2))}{n_{v-1}^2 n_v} \right| \\
&\leq n_k \sum_{v=k+1}^{\infty} \frac{1}{n_{v-1}^2 n_v} \\
&\leq \sum_{v=0}^{\infty} \frac{1}{n_k n_{k+1} 4^v} \\
&= \frac{4}{3 n_k n_{k+1}}.
\end{aligned}$$

Thus, we get

$$\begin{aligned}
|V_{i_1} - V_{i_2}| &\geq \frac{1}{n_k n_{k-1}} - \frac{4}{3 n_k n_{k+1}} \\
&\geq \frac{2}{3 n_k n_{k-1}} \\
&\geq \frac{2}{3 n_{k-1}^{3m+1}} \\
&\geq \frac{2}{3 |i_1 - i_2|^{3m+1}}.
\end{aligned}$$

Therefore, $V$ is a distal sequence with an approximation function

$$Q(x) = \begin{cases} \frac{3n_1^{3m+1}}{2}, & 0 \leq x < n_1; \\ \frac{3x^{3m+1}}{2}, & x \geq n_1. \end{cases}$$

By Theorem 1.22, we have $\operatorname{hull}(V) \cong \operatorname{hull}(L) \cong \Omega$. By Theorem 1.35 there is an $f \in C(\Omega, \mathbb{R})$ such that $(f(T^i(e)))_{i \in \mathbb{Z}} = V$. □

**Remark 2.43.** *For any $r \geq 0$, let*

$$G(x) = \begin{cases} 1, & 0 \leq x < 1; \\ x^r, & x > 1. \end{cases}$$

*It is not hard to see that*

$$h(t) \leq c t^{-4-r}$$

*by choosing* $t_i = t2^{-i-1}$*. The constant* $c$ *depends on* $r$*. It follows that* $G(x)$ *is an approximation function. In particular, this also shows that our* $Q(x)$ *in the above proof is indeed an approximation function.*

### 2.4.3 Pöschel's Results

Here we rewrite some of Pöschel's results from [30], tailored to our purpose.

Let $\tau \geq 1$ be an integer and $\mathfrak{M}$ a Banach algebra of real $\tau$-dimensional sequences $a = (a_i)_{i\in\mathbb{Z}^\tau}$ with the operations of pointwise addition and multiplication of sequences. In particular, the constant sequence 1 is supposed to belong to $\mathfrak{M}$ and have norm one. Moreover, $\mathfrak{M}$ is required to be invariant under translation: if $a \in \mathfrak{M}$, then $\|T_k a\|_{\mathfrak{M}} = \|a\|_{\mathfrak{M}}$ for all $k \in \mathbb{Z}^\tau$, where $T_k a_i = a_{i+k}$.

We denote by $M$ the space of all matrices $A = (a_{i,j})_{i,j\in\mathbb{Z}^\tau}$ satisfying $A_k = (a_{i,i+k}) \in \mathfrak{M}$, $k \in \mathbb{Z}^\tau$, that is, $A_k$ is the $k$-th diagonal of $A$ and it is required to belong to $\mathfrak{M}$. In $M$, we define a Banach space

$$M^s = \{A \in M, \| A \|_s < \infty\}, \quad 0 \leq s \leq \infty,$$

where

$$\|A\|_s = \sup_{k\in\mathbb{Z}^\tau} \|A_k\|_{\mathfrak{M}} e^{|k|s}.$$

Obviously,

$$M^s \subset M^t, \quad \| \cdot \|_s \geq \| \cdot \|_t, \quad 0 \leq t \leq s \leq \infty.$$

In particular, $M^\infty$ is the space of all diagonal matrices in $M$.

**Theorem 2.44.** [30, Theorem A] *Let* $D$ *be a diagonal matrix whose diagonal* $V$ *is a distal sequence for* $\mathfrak{M}$*. Let* $0 < s \leq \infty$ *and* $0 < \sigma \leq \min\{1, \frac{s}{2}\}$*. If* $P \in M^s$ *and* $\| P \|_s \leq \delta \cdot h(\frac{\sigma}{2})^{-1}$*, where* $\delta > 0$ *depends on the dimension* $\tau$ *only, then there exists another diagonal matrix* $\tilde{D}$ *and an invertible matrix* $W$ *such that*

$$W^{-1}(\tilde{D} + P)W = D.$$

*In fact,* $W, W^{-1} \in M^{s-\sigma}$ *and* $\tilde{D} - D \in M^\infty$ *with*

$$\| W - I \|_{s-\sigma}, \| W^{-1} - I \|_{s-\sigma} \leq C \cdot \| P \|_s,$$

$$\| \tilde{D} - D + [P] \|_\infty \leq C^2 \cdot \| P \|_s^2,$$

*where* $C = \delta^{-1} \cdot h(\frac{\sigma}{2})$*, and* $[\cdot]$ *denotes the canonical projection* $M^s \to M^\infty$*. If* $P$ *is self-adjoint, then* $W$ *can be chosen to be orthogonal on* $\ell^2(\mathbb{Z}^\tau)$*. Note that* $h$ *is the function* (2.23) *associated with* $V$*.*

An important consequence of the preceding theorem for discrete Schrödinger operators is the following.

**Theorem 2.45.** [30, Corollary A] *Let $V$ be a distal sequence for some translation invariant Banach algebra $\mathfrak{M}$ of $\tau$-dimensional real sequences. Then for $0 \le \varepsilon \le \varepsilon_0, \varepsilon_0 > 0$ sufficiently small, there exists a sequence $\tilde{V}$ with $\tilde{V} - V \in \mathfrak{M}$, $\| \tilde{V} - V \|_{\mathfrak{M}} \le \frac{\varepsilon^2}{\varepsilon_0^2}$, such that the discrete Schrödinger operator*

$$(\tilde{H}\psi)_i = \varepsilon \sum_{|l|=1} \psi_{i+l} + \tilde{V}_i \psi_i, \quad i \in \mathbb{Z}^\tau$$

*has eigenvalues $\{V_i : i \in \mathbb{Z}^\tau\}$ and a complete set of corresponding exponentially localized eigenvectors with decay rate $1 + \log \frac{\varepsilon_0}{\varepsilon}$.*

Next let us discuss how to apply the above results.

**Pöschel's Example.** Fix $\tau \ge 1$, and let $\mathcal{P}$ be the set of all real $\tau$-dimensional sequences $a = (a_i)$ with period $2^n, n \ge 0$, in each dimension; that is, $a_i = a_j, i - j \in 2^n \mathbb{Z}^\tau$. The closure of $\mathcal{P}$ with respect to the sup norm $\| \cdot \|_\infty$ is a Banach algebra, which we denote by $\mathcal{L}$. It is a subspace of the space of all limit periodic sequences.

Let $\alpha_v, v \ge 1$, be the characteristic function of the set

$$A_v = \begin{cases} \bigcup_{N\in\mathbb{Z}} [N \cdot 2^v, N \cdot 2^v + 2^{v-1}), & v \text{ even;} \\ \bigcup_{N\in\mathbb{Z}} [N \cdot 2^v + 2^{v-1}, N \cdot 2^v + 2^v), & v \text{ odd.} \end{cases}$$

Then, $\alpha_v$ has period $2^v$. Construct an $\tau$-dimensional sequence $V = (V_i)$ such that

$$V_i = \sum_{v=1}^{\infty} \sum_{\mu=1}^{\tau} \alpha_v(i_\mu) 2^{-(v-1)\tau-\mu}, i = (i_1, \cdots, i_\tau) \in \mathbb{Z}^\tau,$$

belongs to $\mathcal{L}$ and lies dense in $[0, 1]$. It is a distal sequence for $\mathcal{L}$ with

$$\| (V - T_k V)^{-1} \|_\infty \le 16^\tau |k|^\tau, \quad 0 \ne k \in \mathbb{Z}^\tau.$$

Applying Theorem 2.45 to this distal sequence $V$, we find that there exists $\tilde{V} \in \mathcal{L}$ and $\varepsilon_0 > 0$ such that for any $0 < \varepsilon \le \varepsilon_0$, the discrete Schrödinger operator with potential $(\frac{\tilde{V}_i}{\varepsilon})_{i\in\mathbb{Z}}$ has the pure point spectrum $\overline{\{\frac{V_i}{\varepsilon} : i \in \mathbb{Z}^\tau\}}$ and a complete set of exponentially localized eigenvectors with decay rate $1 + \log \frac{\varepsilon_0}{\varepsilon}$. Moreover, the spectrum of this Schödinger operator as a set is $\overline{\{\frac{V_i}{\varepsilon} : i \in \mathbb{Z}^\tau\}} = [0, \frac{1}{\varepsilon}]$ since $\overline{\{V_i : i \in \mathbb{Z}^\tau\}} = [0, 1]$.

### 2.4.4 Existence of Uniform Localization

We are now ready to give the proof of Theorem 2.38 about uniform localization. Given a procyclic group $\Omega$ with the condition $\mathscr{A}$, and a minimal translation $T$, we fix a metric $\| \cdot \|$ compatible with the topology. We have already seen that there exists some $f \in C(\Omega, \mathbb{R})$ such that $V = (f(T^i(e)))_{i\in\mathbb{Z}}$ is a distal sequence; compare Lemma 2.42. Clearly, $C(\Omega, \mathbb{R})$ will induce a class of limit-periodic potentials. We denote it by $\mathcal{B}$, and one can check that this class is a translation invariant Banach algebra with the $\ell^\infty$-norm. By Theorem 2.45, there exists a sufficiently small $\varepsilon_0 > 0$ such that for $0 < \varepsilon \leq \varepsilon_0$, there is a sequence $\tilde{V} \in \mathcal{B}$ with $\| \tilde{V} - V \|_\infty \leq \frac{\varepsilon^2}{\varepsilon_0^2}$ so that the discrete Schrödinger operator

$$(H\psi)_i = \psi_{i-1} + \psi_{i+1} + \frac{\tilde{V}_i}{\varepsilon}\psi_i, \quad i \in \mathbb{Z}$$

has eigenvalues $\{\frac{V_i}{\varepsilon}, i \in \mathbb{Z}\}$ and a complete set of corresponding exponentially localized eigenvectors with decay rate $r = 1 + \log \frac{\varepsilon_0}{\varepsilon}$. There exists a sampling function $\tilde{f} \in C(\Omega, \mathbb{R})$ such that $\tilde{f}(T^i(e)) = \frac{\tilde{V}_i}{\varepsilon}$ since $\tilde{V} \in \mathcal{B}$.

For the Schrödinger operator $H$ associated with potential $\tilde{f}(T^i(e))$, denote its matrix representation with respect to the standard orthonormal basis of $\ell^2(\mathbb{Z})$, $\{\delta_n\}_{n\in\mathbb{Z}}$, by the same symbol. Pöschel's theorem also implies that there exists an orthogonal $W : \ell^2(\mathbb{Z}) \to \ell^2(\mathbb{Z})$ (with corresponding matrix denoted by the same symbol) such that

$$H \cdot W = W \cdot D, \tag{2.25}$$

where $D$ is a diagonal matrix with the diagonal $D_0 = (\frac{V_i}{\varepsilon})_{i\in\mathbb{Z}}$. We write $W = (\cdots, W_{-1}, W_0, W_1, \cdots)$ where $W_i$ is the $i$-th diagonal of $W$, and similarly, we write $H = (\cdots, 0, 0, H_{-1}, H_0, H_1, 0, 0, \cdots)$ and $D = (\cdots, 0, 0, D_0, 0, 0, \cdots)$. Moreover, by Theorem 2.44 we have that $W \in M^r$, where $r > 0$ and $M^r$ is a space of matrices associated with the Banach algebra $\mathcal{B}$ (see Subsection 2.4.3 for the description of this space). (Note that $W \in M^r$ follows from [30, Proof of Corollary A].) Since $W \in M^r$, we have $\| W \|_r = \sup_{i\in\mathbb{Z}} \| W_i \|_\infty e^{|i|r} < C$ where $C$ is a constant. So $\| W_i \|_\infty < Ce^{-r|i|}, \forall i \in \mathbb{Z}$. Let $W^{(j)}$ be the $j$-th column of $W$, that is, $W^{(j)}$ is an eigenfunction of $H$. Since $W^{(j)}(k) = W^{(k+(j-k))}(k)$, $W^{(j)}(k)$ is also an entry in $W_{j-k}$, and so $|W^{(j)}(k)| < Ce^{-r|j-k|}$. $C$ is independent of $j$, so the corresponding Schrödinger operator $H$ has ULE. This property is strong enough to imply that the pure point spectrum of $H$ is independent of $\omega$ [25], that is, it is phase stable. In order to see this more explicitly, we would like to prove it in our framework, and furthermore, show that for other $\omega$, the associated Schrödinger operator still has ULE with the same constant $C$. Note that the latter property does not follow from Theorem 2.36.

We need the following lemma.

**Lemma 2.46.** *Suppose we are given matrices $A, B \in \mathbb{R}^{\mathbb{Z}\times\mathbb{Z}}$, one of which has only finitely many non-zero diagonals. Then, we have for the $k$-th diagonal of $Z = AB$,*

$$Z_k = \sum_{l\in\mathbb{Z}} A_l \cdot T^l(B_{k-l}),$$

*where $\cdot$ is the pointwise multiplication (i.e., $A_l \cdot T^l(B_{k-l})$ is still a sequence) and $T$ is the translation defined by $(T(B_{k-l}))_i = (B_{k-l})_{i+1}$ for $i \in \mathbb{Z}$.*

*Proof.* Since for $i, k \in \mathbb{Z}$, we have

$$\begin{aligned} z_{i,i+k} &= \sum_{t\in\mathbb{Z}} a_{i,t} b_{t,i+k} \\ &= \sum_{l\in\mathbb{Z}} a_{i,i+l} b_{i+l,i+k} \\ &= \sum_{l\in\mathbb{Z}} a_{i,i+l} b_{i+l,i+l+k-l}, \end{aligned}$$

the lemma follows. $\square$

Now consider a given $\omega \in \Omega$. By Proposition 1.20 we have $(\tilde{f}(T^i(\omega)))_{i\in\mathbb{Z}} \in \operatorname{hull}((\tilde{f}(T^i(e)))_{i\in\mathbb{Z}})$. If $\omega$ is in the orbit of $e$, that is, $\omega = T^t(e)$ for some $t \in \mathbb{Z}$, ULE with the same constants and eigenvalues follows from unitary operator equivalence directly. However, we write this out in detail so that we see clearly what happens in the case where $\omega$ can only be approximated by elements of the form $T^t(e)$.

By the previous lemma, (2.25) is equivalent to the following form:

$$\forall k \in \mathbb{Z}: \qquad \sum_{l\in\mathbb{Z}} H_l \cdot T^l(W_{k-l}) = \sum_{l\in\mathbb{Z}} W_l \cdot T^l(D_{k-l}).$$

Since $D_j = 0$ for $j \neq 0$ and $H_{\pm 1}$ are both constant equal to one, this simplifies as follows,

$$\forall k \in \mathbb{Z}: \qquad T^{-1} W_{k+1} + H_0 \cdot W_k + T W_{k-1} = W_k \cdot T^k(D_0).$$

If the potential is replaced by $\tilde{f}(T^{i+t}(e))$, with the matrix

$$\tilde{H} = (\cdots, 0, 0, \tilde{H}_{-1}, \tilde{H}_0, \tilde{H}_1, 0, 0, \cdots)$$

such that $\tilde{H}_j(i) = H_j(i+t), j \in \{-1, 0, 1\}$, we still have

$$\forall k \in \mathbb{Z}: \qquad T^{-1}\tilde{W}_{k+1} + \tilde{H}_0 \cdot \tilde{W}_k + T\tilde{W}_{k-1} = \tilde{W}_k \cdot T^k(\tilde{D}_0),$$

where $\tilde{W}_k(i) = W_k(i+t),\ k \in \mathbb{Z}$ and $\tilde{D}_0(i) = D_0(i+t)$. Reversing the steps above, this means that

$$\tilde{H} \cdot \tilde{W} = \tilde{W} \cdot \tilde{D}.$$

We can conclude that $\tilde{H}$ has the pure point spectrum $\overline{\{\frac{V_{i+t}}{\varepsilon} : i \in \mathbb{Z}\}} = \overline{\{\frac{V_i}{\varepsilon} : i \in \mathbb{Z}\}}$. Moreover, $\tilde{W} = (\cdots, \tilde{W}_{-1}, \tilde{W}_0, \tilde{W}_1, \cdots)$ is the eigenfunction matrix of $\tilde{H}$, and for any $i, k \in \mathbb{Z}$, $|\tilde{W}_k(i)| = |W_k(i+t)| \leq Ce^{-r|k|}$. So for the eigenfunction $\tilde{W}^{(j)}$ of $\tilde{H}$, we still have $|\tilde{W}^{(j)}(i)| < Ce^{-r|j-i|}$, and hence ULE with the same constants follows.

If $\lim_{m\to\infty} T^{t_m}(e) = \omega$, that is, $\tilde{f}(T^i(\omega)) = \lim_{m\to\infty} \tilde{f}(T^{i+t_m}(e))$, then for $\tilde{f}(T^{i+t_m}(e))$, we have already seen that

$$\tilde{H}^{(m)} \cdot \tilde{W}^{(m)} = \tilde{W}^{(m)} \cdot \tilde{D}^{(m)}. \tag{2.26}$$

Let $\tilde{W}_k^{(m)}$ be the $k$-th diagonal of $\tilde{W}^{(m)}$, so that $\tilde{W}_k^{(m)}(i) = W_k(i+t_m)$. They are also in the same Banach space, $\mathcal{B}$, which is just a certain class of limit-periodic sequences. So there exists some $\tilde{f}_k \in C(\Omega, \mathbb{R})$ such that $\tilde{W}_k^{(m)}(i) = W_k(i+t_m) = \tilde{f}_k(T^{i+t_m}(e))$. So $\lim_{m\to\infty} \tilde{W}_k^{(m)}(i) = \lim_{m\to\infty} \tilde{f}_k(T^{i+t_m}(e)) = \tilde{f}_k(T^i(\omega))$, and we denote $\tilde{f}_k(T^i(\omega))$ by $\tilde{W}_k^{(\infty)}(i)$. Similarly, $\lim_{m\to\infty} \tilde{D}^{(m)}$ exists and $\tilde{D}_0^{(\infty)}(i) = f(T^i(\omega))$, where $\tilde{D}_0^{(\infty)}$ is the 0-th diagonal of $\tilde{D}^{(\infty)}$. Thus, as we let $m \to \infty$, (2.26) takes the following form:

$$\tilde{H}^{(\infty)} \cdot \tilde{W}^{(\infty)} = \tilde{W}^{(\infty)} \cdot \tilde{D}^{(\infty)}, \tag{2.27}$$

where $\tilde{H}^{(\infty)}$ is (the matrix representation of) the Schrödinger operator with potential $\tilde{f}(T^i(\omega))$. Equation (2.27) implies that $\tilde{H}^{(\infty)}$ has the pure point spectrum $\overline{\{\frac{V_i}{\varepsilon} : i \in \mathbb{Z}\}}$, and its eigenfunctions are uniformly localized since $|(\tilde{W}^{(\infty)})^{(j)}(k)| < Ce^{-r|j-k|}$ for any $j, k \in \mathbb{Z}$, where $(\tilde{W}^{(\infty)})^{(j)}$ is the $j$-th column of $\tilde{W}^{(\infty)}$. This completes the proof of Theorem 2.38.

# Chapter 3

# Open Problems

We conclude our thesis with a number of interesting open problems concerning the spectral properties of limit-periodic Schrödinger operators.

Given a minimal translation $T$ of a procyclic group $\Omega$, consider for $f \in C(\Omega, \mathbb{R})$ and $\omega \in \Omega$ the spectral type of the associated Schrödinger operator $H_\omega$ with potential given by $V_\omega(n) = f(T^n(\omega))$.

***Problem 1.*** Is it true that for $f$ from a suitable dense subset of $C(\Omega, \mathbb{R})$, $H_\omega$ has pure point spectrum for (Haar-) almost every $\omega \in \Omega$?

We already know that for generic $f \in C(\Omega, \mathbb{R})$, $H_\omega$ has purely singular continuous spectrum for every $\omega \in \Omega$, and also that for $f$ from a suitable dense subset of $C(\Omega, \mathbb{R})$, $H_\omega$ has purely absolutely continuous spectrum for every $\omega \in \Omega$. Thus, an affirmative answer to Problem 1 would clarify the effect of the choice of $f$ on the spectral type. Since the methods of Pöschel are essentially restricted to large potentials, one should not expect them to yield an answer to Problem 1 and one should in fact pursue methods involving some randomness aspect.

Note, however, the different quantifier on $\omega$ in Problem 1, compared to the results just quoted. We exhibit $(\Omega, T, f)$ for which $H_\omega$ has pure point spectrum for every $\omega \in \Omega$. From this perspective, the following problem arises naturally:

***Problem 2.*** Is the spectral type of $H_\omega$ always the same for every $\omega \in \Omega$?

For quasi-periodic potentials, this is known not to be the case (cf. [26]). However, the mutual approximation by translates for two given elements in the hull is stronger in the limit-periodic case than in the quasi-periodic case, so it is not clear if similar counterexamples to uniform spectral types exist in the limit-periodic world.

Another related problem is the following:

***Problem 3.*** Is the spectral type of $H_\omega$ always pure?

Again, in the quasi-periodic world, this is known not to be the case: there are examples that have both absolutely continuous spectrum and point spectrum (cf. [4, 5]).

Returning to the issue of point spectrum, one interesting aspect of the result stated (in the continuum case) by Molchanov and Chulaevsky in [29] is the coexistence of pure point spectrum with the absence of non-uniform hyperbolicity. That is, in their examples, the Lyapunov exponent vanishes on the spectrum and yet the spectral measures are pure point. This is the only known example of this kind and it would therefore be of interest to have a complete published proof of a result exhibiting this phenomenon. Especially since our study is carried out in a different framework, we ask within this framework the following question:

***Problem 4.*** For how many $f \in C(\Omega, \mathbb{R})$ does the Lyapunov exponent vanish throughout the spectrum and yet $H_\omega$ has pure point spectrum for (almost) every $\omega \in \Omega$?

Given the existing ideas, it is conceivable that Problems 1 and 4 are closely related and may be answered by the same construction. If this is the case, it will then still be of interest to show for a dense set of $f$'s that there is almost sure pure point spectrum with *positive* Lyapunov exponents.